\documentclass[11pt,leqno]{amsart}
\usepackage{a4wide}
\usepackage{verbatim,amssymb,amsfonts,latexsym,amsmath,amsthm,bbm,nicefrac,xspace,enumerate,tikz}
\usepackage[normalem]{ulem}

\newtheorem{theorem}{Theorem}[section]

\newtheorem{lemma}[theorem]{Lemma}
\newtheorem{corollary}[theorem]{Corollary}
\newtheorem{proposition}[theorem]{Proposition}
\newtheorem{question}[theorem]{Question}
\newtheorem{claim}[theorem]{Claim}
\theoremstyle{definition}

\newtheorem{definition}[theorem]{Definition}

\newtheorem{remark}[theorem]{Remark}
\newcommand{\Q}{\mathbb{Q}}
\newcommand{\R}{\mathbb{R}}
\newcommand{\IR}{\mathbb{R}}

\newcommand{\IN}{\mathbb{N}}
\newcommand{\IQ}{\mathbb{Q}}
\newcommand{\IZ}{\mathbb{Z}}

\newcommand{\CC}{\mathbb{C}}

\newcommand{\F}{\mathcal{F}}

\newcommand{\K}{\mathcal{K}}

\newcommand{\U}{\mathcal{U}}
\newcommand{\V}{\mathcal{V}}

\newcommand{\Ra}{\Rightarrow}
\newcommand{\defeq}{\overset{\mbox{\tiny\sf def}}=}
\newcommand{\explicitSet}[1]{\left\lbrace #1 \right\rbrace}

\newcommand{\set}[2]{\explicitSet{#1 \colon #2}}

\renewcommand{\a}{\alpha}
\renewcommand{\b}{\beta}

\newcommand{\e}{\varepsilon}

\renewcommand{\k}{\kappa}
\newcommand{\s}{\sigma}

\newcommand{\w}{\omega}

\newcommand{\sub}{\subseteq}

\newcommand{\cl}{\mathrm{cl}}
\newcommand{\closure}[1]{\overline{#1}}

\newcommand{\cf}{\mathsf{cf}}
\newcommand{\card}[1]{\left\lvert #1 \right\rvert}

\newcommand{\PP}{\mathbb{P}}

\newcommand{\BB}{\mathbb{B}}

\newcommand{\continuum}{\mathfrak{c}}

\newcommand{\dom}{\mathfrak d}
\newcommand{\bdd}{\mathfrak b}

\newcommand{\scr}[1]{\ensuremath{\mathcal {#1}}}

\newcommand{\cov}[1]{\mathsf{cov}(\scr{#1})}

\newcommand{\ch}{\ensuremath{\mathsf{CH}}\xspace}

\newcommand{\zfc}{\ensuremath{\mathsf{ZFC}}\xspace}

\newcommand{\ma}{\ensuremath{\mathsf{MA}}\xspace}
\newcommand{\masc}{\ensuremath{\mathsf{MA}(\s\text{-centered})}\xspace}

\begin{document}

\title{First-countable Lindel\"{o}f scattered spaces}
\author{Taras Banakh}
\address {
T. O. Banakh\\
Ivan Franko National University of Lviv (Ukraine) and Jan Kochanowski University in Kielce (Poland)}
\email{t.o.banakh@gmail.com}
\author{Will Brian}
\address {
W. R. Brian\\
Department of Mathematics and Statistics\\
University of North Carolina at Charlotte\\
Charlotte, NC 
(USA)}
\email{wbrian.math@gmail.com}
\urladdr{wrbrian.wordpress.com}
\author{Alejandro R\'ios Herrej\'on}
\address {
A. R. Herrej\'{o}n\\
Departamento de Matem\'aticas, Facultad de Ciencias, Universidad Nacional Aut\'onoma de M\'exico, Circuito ext. s/n, Ciudad Universitaria, C.P. 04510,  M\'exico, CDMX}
\email{chanchito@ciencias.unam.mx}
\thanks{The third author was supported by CONACYT grant no. 814282.}




\begin{abstract}
We study the class of first-countable Lindel\"of scattered spaces, or ``FLS'' spaces. 
While every $T_3$ FLS space is homeomorphic to a scattered subspace of $\mathbb Q$, the class of $T_2$ FLS spaces turns out to be surprisingly rich. 
Our investigation of these spaces reveals close ties to $Q$-sets, Lusin sets, and their relatives, and to the cardinals $\mathfrak{b}$ and $\mathfrak{d}$. Many natural questions about FLS spaces turn out to be independent of $\mathsf{ZFC}$.

We prove that there exist uncountable FLS spaces with scattered height $\omega$. On the other hand, an uncountable FLS space with finite scattered height exists if and only if $\mathfrak{b} = \aleph_1$. 
We prove some independence results concerning the possible cardinalities of FLS spaces, and concerning what ordinals can be the scattered height of an FLS space.
Several open problems are included.
\end{abstract}

\maketitle

\section{Introduction}

In this paper we introduce and study the class of first-countable Lindel\"{o}f scattered spaces, henceforth referred to as \emph{FLS} spaces. Here, and throughout, a \emph{space} means a Hausdorff topological space, unless otherwise specified.

Our starting point is a result of Marlene Gewand from 1984:
\begin{theorem}[Gewand \cite{Gewand}]
Every regular FLS space is countable.
\end{theorem}
\begin{corollary}\label{cor:Gewand}
A topological space is a regular FLS space if and only if it is homeomorphic to a scattered subspace of $\Q$.
\end{corollary}
\begin{proof}
To prove the nontrivial direction, suppose $X$ is a regular FLS space. Gewand's theorem states that such a space is countable. Clearly, a space that is both countable and first-countable is second-countable. Second-countable regular spaces are metrizable, by the Urysohn Metrization Theorem. And finally, every countable metrizable space is homeomorphic to a subspace of $\Q$ (by a theorem of Fr\'echet; see \cite[Exercise 4.3.H(b)]{Engelking}).
\end{proof}

Thus it seems that the regular FLS spaces are rather tame and simple, and easy to understand. 
This raises a question: what about non-regular FLS spaces?
As we hope to demonstrate in what follows, the non-regular FLS spaces form an unexpectedly rich class, where many natural questions have surprising answers, some independent of \zfc.
Our investigation of these spaces reveals their close ties to $Q$-sets, Lusin sets and their relatives, and to the cardinal characteristics $\mathfrak{b}$ and $\mathfrak{d}$.

Throughout the paper, all our results are aimed at exploring two basic questions:
\begin{itemize}
\item[\small{$(\mathrm{Q}1)$}] What are the possible cardinalities of FLS spaces?
\item[\small{$(\mathrm{Q}2)$}] What ordinals are the scattered height of an FLS space?
\end{itemize}

Concerning the first question {\small{$(\mathrm{Q}1)$}}, the main results include:
\begin{itemize}
\item[$\bullet$] There are FLS spaces of cardinality $\mathfrak{c}$. (Theorem~\ref{thm:Existence})
\item[$\bullet$] It is consistent with arbitrarily large values of $\mathfrak{c}$ that there are FLS spaces of cardinality $\k$ for every $\k \leq \mathfrak c$. (Corollary~\ref{c:size-d1})
\end{itemize}
The cardinality $\continuum$ in these results is as large as possible: a celebrated theorem of Arhangel'ski\u{\i} from \cite{Arhangelskii} states that every first-countable Lindel\"{o}f space (hence every FLS space) has cardinality at most $\continuum$. 

Restricting our attention to FLS spaces with finite scattered height, the cardinality question becomes more subtle. We show that:
\begin{itemize}
\item[$\bullet$] There is an uncountable FLS space of finite scattered height if and only if $\bdd = \aleph_1$. (Corollary~\ref{c:b})
\item[$\bullet$] Every FLS space with finite scattered height has cardinality $\leq\!\dom$. (Corollary~\ref{c:bd})
\item[$\bullet$] It is consistent with arbitrary values of $\dom$ and $\continuum$ that there are FLS spaces with finite scattered height of every cardinality $\leq\!\dom$. (Theorem~\ref{thm:CohenModel})
\item[$\bullet$] It is consistent with arbitrarily large values of $\dom$ that there are FLS spaces with finite scattered height of cardinality $\aleph_1$, but no larger. (Theorem~\ref{thm:DualCohen})
\end{itemize}
Our main unanswered question regarding the cardinalities of FLS spaces is whether it is consistent to have some cardinal $\k < \mathfrak{c}$ that is not the cardinality of an FLS space.

\vspace{2mm}

Concerning the second question {\small{$(\mathrm{Q}2)$}} listed above, we prove: 
\begin{itemize}
\item[$\bullet$] If $\a$ is the scattered height of an FLS space, then $\a < \continuum^+$ and $\cf(\a) \leq \w$. (Proposition~\ref{p:ImpossibleRanks})
\item[$\bullet$] It is consistent with arbitrary values of $\mathfrak{c}$ that there is an FLS space of scattered height $\a$ if and only if $\a < \continuum^+$ and $\cf(\a) \leq \w$. (Corollary~\ref{c:height})
\end{itemize}
Note that the restrictions listed in the first bullet point are easy to prove, given that every FLS space is Lindel\"{o}f and has cardinality $\leq\mathfrak{c}$.
The model witnessing the second bullet point is obtained simply by adding any number of Cohen reals to a model of \ch. 
Thus, in other words, the Cohen model contains FLS spaces of every scattered height not obviously impossible for a Lindel\"{o}f space of cardinality $\leq\!\continuum$.
It is a sort of ``maximal'' model for the question of what ordinals are the scattered height of an FLS space.
Our main open question in this direction is whether there is also a ``minimal'' model, where an ordinal is the scattered height of an FLS space only if some \zfc theorem requires it to be such. Our main progress in this direction is:
\begin{itemize}
\item[$\bullet$] Martin's Axiom implies that $\a < \continuum^+$ is the scattered height of an FLS space if and only if either $\cf(\a) \leq \w$, or else $\a = \lambda+n$ for some $n \in \w$ and some $\lambda$ with $\cf(\lambda) = \mathfrak{c}$. (Corollary~\ref{c:height2})
\end{itemize}
It is necessary for ordinals of the first kind to be the scattered height of an FLS space (Theorem~\ref{t:construct3}), but we do not know whether it is necessary for ordinals of the second kind (Question~\ref{q:possibleheights}).

Some of the results listed here in this brief summary have been simplified, and a fuller (but more complicated) statement of these results can be found below. There are also several results and questions found in what follows that are not mentioned here in the introduction.

\vspace{2mm}

Recall that a space is \emph{hereditarily Lindel\"{o}f} if each of its subspaces is Lindel\"{o}f. It is not difficult to prove that every hereditarily Lindel\"{o}f scattered space is countable. Thus none of the more interesting FLS spaces we investigate here are hereditarily Lindel\"{o}f. However, we will repeatedly make use of the fact that every closed subspace of a Lindel\"{o}f space is Lindel\"{o}f. Because every subspace of a first-countable space is first-countable and every subspace of a scattered space is scattered (i.e., these are hereditary properties), it follows that every closed subspace of an FLS space is itself an FLS space.

\vspace{2mm}

The organization of the paper is as follows. Section 2 contains some preliminary matter, mostly to set our notation and conventions, but also to collect some known results into one place for the reader's convenience. Section 3 contains a short development of our general method for constructing FLS spaces: scattered modifications of topological spaces. Section 4 contains the relevant applications of this method, in the form of four different constructions of FLS spaces. 
In Section 5, we prove several theorems concerning general properties of FLS spaces. This section contains many results constraining the possible combinations of cardinalities and scattered heights of FLS spaces. Finally, Section 6 contains both a synthesis of the results proved in Sections 4 and 5, and a few forcing constructions showing some ways in which these results cannot be improved.

\section{Preliminaries}

\subsection{Some notation} Let $\IR_+$ denote the set $[0,\infty)$ of non-negative real numbers.

For a function $f:X\to Y$ between sets, element $x\in X$ and subset $A\subseteq X$ we denote by $f(x)$ the value of $f$ at $x$ and by $f[A]$ the image $\{f(x):x\in A\}$ of the set $A$ under the map $f$. Distinguishing between $f(x)$ and $f[A]$ is important for functions on ordinals (which are simultaneously sets of ordinals). Also for an element $y\in Y$ and subset $B\subseteq Y$ we put $f^{-1}\{y\}\defeq\{x\in X:f(x)=y\}$ and $f^{-1}[B]\defeq\{x\in X:f(x)\in B\}$.

For an element $x$ of a partially ordered set $(X,\le)$ we define the lower and upper sets
$${\downarrow}x\defeq\{y\in X:y\le x\}\quad\mbox{and}\quad{\uparrow}x\defeq\{y\in X:x\le y\}$$
of $x$, respectively.

\subsection{Ordinals} Every ordinal is identified with the set of smaller ordinals. An ordinal $\alpha$ is 
\begin{itemize}
\item {\em successor} if $\alpha=\beta+1\defeq\beta\cup\{\beta\}$ for some ordinal $\beta$;
\item {\em limit} if $\alpha$ is not successor.
\end{itemize}
According to this definition, the least ordinal $0$ is limit. By $\w$ be denote the smallest nonzero limit ordinal. Elements of $\w$ are called {\em finite ordinals} or {\em natural numbers}. Let $\IN\defeq\w\setminus\{0\}$ denote the set of positive integers. 

Every ordinal $\alpha$ can be uniquely written as the sum $\lambda+n$ of a limit ordinal $\lambda$ and a finite ordinal $n$.
This limit ordinal $\lambda$ is called the {\em limit part} of $\alpha$ and denoted $\lim\alpha$. The limit part of any finite ordinal is zero.
The {\em cofinality} of an ordinal $\alpha$ is denoted $\cf(\alpha)$; let us note that if $\alpha$ is successor then $\cf(\alpha)=1$.

For a set $X$ its {\em cardinality} $|X|$ is defined as the smallest ordinal $\alpha$ that admits a bijective function $f:\alpha\to X$. 
A set $X$ is {\em countable} if $|X| \le \w$. So, countable sets are either finite or countably infinite. For a cardinal $\kappa$ we denote by $\kappa^+$ the successor cardinal of $\kappa$.

\subsection{Scattered spaces}

Given a topological space $X$, the \emph{Cantor-Bendixson derivative} of $X$, denoted $X^{(1)}$, is the set of all non-isolated points in $X$.
Iterating the Cantor-Bendixson derivative, we obtain by transfinite recursion an ordinal-indexed sequence of subspaces of $X$: we define $X^{(0)} = X$ and for all ordinals $\a > 0$,
$$
X^{(\a)} = \textstyle \bigcap_{\xi < \a} \left( X^{(\xi)} \right){\!}^{(1)}.$$
A straightforward transfinite induction shows that these $X^{(\a)}$'s form a decreasing sequence of closed subsets of $X$. Let $X^{(\infty)}=\bigcap_{\alpha}X^{(\alpha)}$ and $\hbar:X\to \hbar[X]$ be the function assigning to every $x\in X\setminus X^{(\infty)}$ the unique ordinal $\alpha$ such that $x\in X^{(\alpha)}\setminus X^{(\alpha+1)}$, and to every $x\in X^{(\infty)}$ the smallest ordinal $\alpha$ such that $X^{(\alpha)}=X^{(\infty)}$. The ordinal $\hbar(x)$ is the {\em Cantor-Bendixson height} of the point $x$ in $X$. The function $\hbar:X\to \hbar[X]$ is the {\em height function} for $X$ and its range $\hbar[X]$ is an ordinal, the {\em Cantor--Bendixson height} of the space $X$.

A topological space $X$ is {\em scattered} if every nonempty subspace of $X$ has an isolated point. It is well known (and easy to prove) that $X$ is scattered if and only if $X^{(\infty)}=\emptyset$ if and only if $X^{(\hbar[X])}=\emptyset$. In this case, $\hbar[X]$ is also called the \emph{scattered height} of $X$. 

A topological space $X$ is called {\em crowded} if it has no isolated points. For every topological space $X$ the set $X^{(\infty)}$ is the largest crowded subspace of $X$.

\subsection{Countability properties of topological spaces}

Recall that a space $X$ is
\begin{itemize}
\item {\em Lindel\"of} if every open cover of $X$ has a countable subcover;
\item {\em hereditarily Lindel\"of} if each subspace of $X$ is Lindel\"of;
\item {\em separable} if $X$ contains a countable dense subset;
\item {\em hereditarily separable} if each subspace of $X$ is separable;
\item {{\em perfect} if every closed subset is of type $G_\delta$ in $X$;}
\item of {\em countable pseudocharacter} {if every singleton $\{x\}$ is of type $G_\delta$ in $X$;}
\item of {\em countable extent} if every closed discrete subspace of $X$ is countable;
\item {of {\em countable separating weight} if there exists a countable family $\U$ of open sets in $X$ such that for any distinct points $x,y\in X$ some set $U\in\U$ contains $x$ but not $y$.}
\end{itemize}

\subsection{Generalized metric spaces} In this subsection we collect the necessary information on generalized metric spaces. For more details, see \cite[\S9,10]{Grue}.

Let $X$ be a set. For a function $d:X\times X\to \IR_+$ consider the following conditions, where $x,y,z$ are arbitrary elements of $X$:
\begin{itemize}
\item[(i)] $d(x,x)=0$;
\item[(ii)] $d(x,y)=0$ if and only if $x=y$; 
\item[(iii)] $d(x,y)=d(y,x)$;
\item[(iv)] $d(x,z)\le d(x,y)+d(y,z)$.
\end{itemize}
A function $d:X\times X\to\IR_+$ is called
\begin{itemize}
\item a {\em premetric} if it satisfies the condition (i);
\item a {\em symmetric} if it satisfies the conditions (ii) and (iii);
\item a {\em quasi-metric} if it satisfies the conditions (ii) and (iv);
\item a {\em metric} if it is a symmetric and a quasi-metric.
\end{itemize}

We say that the topology $\tau$ of a topological space $X$ is generated by a premetric $d$ if $\tau$ coincides with the family of all sets $U\subseteq X$ such that for every $x\in U$ there exists a positive real $\e$ such that $B_d(x;\e)\defeq\{y\in X:d(x,y)<\e\}\subseteq U$. 

A topological space is called {\em metrizable} (resp. {\em symmetrizable, quasi-metrizable\/}) if its topology is generated by some metric (resp. symmetric, quasi-metric). It is easy to see that every quasi-metrizable space $X$ is first-countable, as for every $x\in X$ the countable family $\{B_d(x;\frac1n):n\in\IN\}$ is a neighborhood base at $x$. {It is known \cite[\S9]{Grue} that every symmetrizable first-countable  Hausdorff space is perfect.}

A topological space $X$ is {\em submetrizable} if there is a continuous metric $d: X \times X \to \R_+$ (though this metric does not need to generate the topology).
Every submetrizable space $X$ admits a continuous bijective map onto a metrizable space (namely, $X$ with the topology generated by the metric $d$ in the definition of submetrizability). It follows that every submetrizable space is {\em functionally Hausdorff\/}, meaning that for every distinct points $x,y\in X$ there exists a continuous function $f:X\to\IR$ such that $f(x)\ne f(y)$.

A subset $U$ in a topological space $X$ is called {\em regular open} if $U$ coincides with the interior of its closure. A topological space $X$ is called 
\begin{itemize}
\item {\em regular} if any neighborhood of any point $x\in X$ contains a closed neighborhood of $x$; 
\item {\em almost regular} if any regular open neighborhood of any point $x\in X$ contains a closed neighborhood of $x$; 
\item {\em semiregular} if regular open sets form a base of the topology of $X$.
\end{itemize}
It is easy to see that a topological space is regular if and only if it is both almost regular and semiregular.

For a topological space $X$ its {\em semiregularization} $X_{sr}$ is the space $X$ endowed with the topology generated by the base consisting of all regular open sets. It is known \cite[1.7.8]{Engelking} that the semiregularization of any topological space is semiregular. A topological space is called {\em nearly metrizable} if its semiregularization is metrizable, see \cite{MM}. The continuity of the identity map $X\to X_{sr}$ implies that nearly metrizable spaces are submetrizable.

For a subset $M$ of a topological space $(X,\tau)$, let $\tau_M$ be the topology on $X$ generated by the base $\tau\cup\{\{a\}:a\in X\setminus M\}$. The topological space $X_M= (X,\tau_M)$ is called the {\em Bing--Hanner modification} of $X$ and is well-known as one of standard counterexamples in General Topology, see \cite[5.1.22]{Engelking} and \cite[Example 71]{SS}. The following metrizability criterion for the Bing--Hanner modifications should be known but we could not find a corresponding reference. 

\begin{lemma}\label{l:Michael} Let $X$ be a metrizable topological space and $M$ be a subset of $X$. The space $X_M$ is metrizable if and only if $M$ is a $G_\delta$-set in $X$.
\end{lemma}

\begin{proof} First assume that $M$ is a $G_\delta$-set in $X$ and write the complement $X\setminus M$ as the union $\bigcup_{n\in\w}F_n$ of an increasing sequence $(F_n)_{n\in\w}$ of closed sets in $X$ such that $F_0=\emptyset$. Let $Y\defeq X/M=\{M\}\cup(X\setminus M)$ be the quotient set of $X$ by its subset $M$, and $q:X\to Y$ be the quotient function. Let  $\lambda:Y\to [0,1]$ be a unique function such that $\lambda(M)=0$ and $\lambda^{-1}\{2^{-n}\}=F_{n+1}\setminus F_n$ for every $n\in\w$. The function $\lambda$ determines the metric $d:Y\times Y\to\IR_+$,
$$d:(x,y)\mapsto\begin{cases}
\max\{\lambda (x),\lambda (y)\}&\mbox{if $x\ne y$};\\
0&\mbox{otherwise};
\end{cases}
$$ 
on the set $Y$. It is easy to check that the function $g:X_M\to X\times Y$, $g:x\mapsto(x,q(x))$, is a topological embedding of the space $X_M$ into the metrizable space $X\times Y$, witnessing that the space $X_M$ is metrizable.
\smallskip

Now assume conversely that the space $X_M$ is metrizable. Then the open subspace $X\setminus M$ of $X_M$ can be written as the union $\bigcup_{n\in\w}F_n$ of a sequence of closed subsets in $X_M$. The definition of the topology of $X_M$ guarantees that for every $n\in\w$ the closure $\bar F_n$ of $F_n$ in $X$ does not intersect the set $M$, which implies that $X\setminus M=\bigcup_{n\in\w}\bar F_n$ is an $F_\sigma$-set in $X$ and hence $M$ is a $G_\delta$-set in $X$.
\end{proof}

\subsection{Some cardinal characteristics of the continuum}

If $f,g \in \w^\w$, we write $f \leq g$ to mean that $f(n) \leq g(n)$ for all $n \in \w$, and we write $f \leq^* g$ to mean that $f(n) \leq g(n)$ for all but finitely many $n \in \w$. A subset $B$ of $\w^\w$ is $\leq^*$-\emph{bounded} if there is some $g \in \w^\w$ such that $f \leq^* g$ for all $f \in B$. The notion of a $\leq$-bounded subset of $B$ is defined similarly. While there are countable $\leq$-unbounded subsets of $\w^\w$ (e.g., the set of all constant functions), a straightforward diagonalization shows that every $\leq^*$-unbounded subset of $\w^\w$ is uncountable. The \emph{unbounding number} $\bdd$ is defined to be the least cardinality of a $\leq^*$-unbounded subset of $\w^\w$. The {\em dominating number} $\mathfrak d$ is the smallest cardinality of a subset $D\subseteq\w^\w$ such that for every $f\in\w^\w$ there exists $g\in\mathcal D$ with $f\le^* g$.

For two functions $f,g\in\w^\w$ we write $f=^*g$ if the set $\{n\in\w:f(n)\ne g(n)\}$ is finite.

\begin{lemma}\label{lem:RemoveTheStar} A subset $B\subseteq\w^\w$ is $\le^*$-bounded if and only if there exists a countable set $C\subseteq \w^\w$ such that for every $x\in B$ there exists $y\in C$ with $x\le y$.
\end{lemma}

\begin{proof} If $B$ is $\le^*$-bounded, then there exists $y\in\w^\w$ such that $x\le^* y$ for every $x\in B$. It is easy to see that the set $C{\defeq}\{x\in \w^\w:x=^* y\}$ is countable and for every $x\in B$ the function $z=\max\{x,y\}$ belongs to the set $C$ and has $x\le z$.

Now assume that $C$ is a countable set such that $\forall x\in B\;\exists y\in C\;(x\le y)$. Since countable sets in $\w^\w$ are $\le^*$-bounded, there exists $y\in\w^\w$ such that $x\le^* y$ for every $x\in C$. Observe that for every $x\in B$ there exists $z\in C$ with $x\le z$. Then $x\le z\le^* y$ implies $x\le^* y$, which means that the set $B$ is $\le^*$-bounded.
\end{proof} 


\begin{lemma}\label{l:d} Let $\F$ be a family of compact subsets in $\w^\w$. If $|\F|<\mathfrak d$, then the union $\bigcup\F$ is disjoint with some uncountable compact subset of $\w^\w$.
\end{lemma}

\begin{proof} For every compact set $F\in\F$, find a function $f_F\in\w^\w$ such that $F\subseteq{\downarrow}f_F$. By the definition of the cardinal $\mathfrak d>|\F|$, there exists a function $g\in\w^\w$ such that $g\not\le^*\! f_F$ for every $F\in\F$. Then $K\defeq \{x\in \w^\w:g\le x\le g+1\}$ is an uncountable compact subset, disjoint with $\bigcup_{F\in\F}{\downarrow}f_F\supseteq\bigcup\F$.
\end{proof}

By $\cov{M}$ we denote the smallest cardinality of a cover of $\w^\w$ by nowhere dense subsets. It is clear that $\cov{M}\le\mathfrak d$. On the other hand, the strict inequality $\cov{M}<\mathfrak d$ holds in some models of \zfc, see \cite[\S11]{Blass}. By Theorem 2.4 in \cite{Blass},
$$\cf(\mathfrak b)=\mathfrak b\le\cf(\mathfrak d)\le\mathfrak d.$$

A topological space $X$ is called a {\em $Q$-space} if every subset in $X$ is of type $F_\sigma$. Let $\mathfrak q_0$ be the smallest cardinality of a metrizable separable space which is not a $Q$-space. It is known  that
$$\w_1\le\mathfrak p\le\mathfrak q_0\le\mathfrak b,$$
where $\mathfrak p$, the \emph{pseudo-intersection number}, denotes the smallest cardinality of a family $\U$ of infinite subsets of $\w$ such that for any finite subfamily $\F\subseteq\U$ the intersection $\bigcap\F$ is infinite but $\U$ has no infinite pseudointersection in the sense that for every infinite set $I\subseteq \w$ there exists a set $U\in\U$ such that $I\setminus U$ is infinite.
Martin's Axiom implies $\mathfrak p=\mathfrak c$. For more information on $Q$-sets and the cardinal $\mathfrak q_0$, see \cite{FM}, \cite{Miller} or \cite{BMZ}.

Some ``non-metrizable'' versions of the cardinal $\mathfrak q_0$ were considered in \cite{BB}. For $i\in\{1,2,3\}$ let $\mathfrak q_i$ be the smallest cardinality of a second-countable $T_i$-space, which is not a $Q$-space. Since second-countable $T_3$-spaces are metrizable, the cardinal $\mathfrak q_3$ coincides with $\mathfrak q_0$.  By \cite{BB}, we have the inequalities
$$\mathfrak p\le\mathfrak q_1\le\mathfrak q_2\le\mathfrak q_3=\mathfrak q_0.$$



\subsection{Special subsets in topological spaces} In this section we recall the necessary information on special sets. More details can be found in \cite{Miller}.

A subset $B$ of a topological space $X$ is called {\em Bernstein} if for every uncountable compact set $K\subseteq B$ the sets $K\cap B$ and $K\setminus B$ are not empty.

\begin{lemma}\label{lem:Bernstein}
For every cardinal $\k$ with $2 \leq \k \leq \continuum$, there is a partition of the real line into $\k$ Bernstein sets.
\end{lemma}

\begin{proof}
This well-known result, originally due to Felix Bernstein, can be proved by transfinite induction. A proof is sketched in chapter 4 of \cite{Oxtoby}.
\end{proof}

A subset $L$ of a topological space $X$ is called
\begin{itemize}
\item {\em Lusin} if for any nowhere dense set $A\subseteq X$ the intersection $A\cap L$ is countable;
\item {\em $\K$-Lusin} if for any compact set $K\subseteq X$ the intersection $K\cap L$ is countable;
\item {\em generalized $\K$-Lusin} if  $|K\cap L|<|L|$ for any compact set $K\subseteq X$;
\item {\em generalized $\K_\sigma$-Lusin} if  $|S\cap L|<|L|$ for any $\sigma$-compact set $S\subseteq X$.
\end{itemize}
From now on, by a ({\em generalized\/}) {\em $\K$-Lusin set} we mean a (generalized) $\K$-Lusin set in $\w^\w$.

Observe that each Lusin set in $\w^\w$ is $\K$-Lusin and each uncountable $\K$-Lusin set is generalized $\K$-Lusin. On the other hand, each generalized $\K$-Lusin set of cardinality $\aleph_1$ is $\K$-Lusin. Every generalized $\K_\sigma$-Lusin set is generalized $\K$-Lusin and every generalized $\K$-Lusin set of cardinality $\kappa$ with $\cf(\kappa)>\w$ is generalized $\K_\sigma$-Lusin.

Since generalized $\K$-Lusin sets will play an important role in evaluation of possible cardinalities and scattered heights of FLS spaces, let us prove two  elementary (but helpful) results on possible cardinalities of generalized $\K$-Lusin sets. 

\begin{lemma}\label{l:genL} There exists a generalized $\K$-Lusin set of cardinality $\mathfrak b$.
\end{lemma}

\begin{proof} The definition of the cardinal $\mathfrak b$ yields a $\le^*$-unbounded set $\{f_\alpha\}_{\alpha\in\mathfrak b}$ in $\w^\w$. Since every set of cardinality $<\mathfrak b$ in $\w^\w$ is $\le^*$-bounded, we can construct inductively a transfinite sequence $(g_\alpha)_{\alpha\in\mathfrak b}$ in $\w^\w$ satisfying the following conditions  for every $\alpha\in\mathfrak b$:
\begin{itemize}
\item $\forall \beta\in\alpha\;\;(g_\beta\le^* g_\alpha)$;
\item $f_\alpha\le g_\alpha$.
\end{itemize} The last condition ensures that the set $L=\{g_\alpha\}_{\alpha\in\mathfrak b}$ is $\le^*$-unbounded and hence has cardinality $\mathfrak b$. We claim that this set is generalized $\K$-Lusin.  Indeed, for any compact set $K\subseteq \w^\w$ we can find a function $g\in\w^\w$ such that $K\subseteq{\downarrow} g$. Since the set $\{f_\alpha\}_{\alpha\in\mathfrak b}$ is $\le^*$-unbounded, there exists $\alpha\in\mathfrak b$ such that $f_\alpha\not\le^* g$. We claim that $L\cap K\subseteq\{g_\beta:\beta<\alpha\}$. Assuming that $g_\beta\in L\cap K\subseteq L\cap {\downarrow}g $ for some $\beta\ge\alpha$, we conclude that $f_\alpha\le g_\alpha\le^* g_\beta\le g$ and hence $f_\alpha\le^* g$, which contradicts the choice of $\alpha$. Thefore, the intersection $K\cap L\subseteq\{g_\beta:\beta<\alpha\}$ has cardinality $\le|\alpha|<\mathfrak b=|L|$, witnessing that the set $L$ is generalized $\K$-Lusin.
\end{proof}



\begin{lemma}\label{l:Kbd} Let $L\subseteq\w^\w$ be a generalized $\K$-Lusin set of cardinality $\kappa$.
\begin{enumerate}
\item[\textup{(1)}] If $\cf(\kappa)>\w$, then $\kappa\ge\mathfrak b$;
\item[\textup{(2)}] $\cf(\kappa)\le\mathfrak d$.
\end{enumerate}
\end{lemma}

\begin{proof} (1) To derive a contradiction, assume that $\cf(\kappa)>\w$ and $|L|=\kappa<\mathfrak b$. By the definition of the cardinal $\mathfrak b$, the set $L$ is $\le^*$-bounded and by Lemma~\ref{lem:RemoveTheStar}, $L\subseteq\bigcup_{f\in F}{\downarrow}f$ for some countable set $F$. Since $\cf(\kappa)>\w$, for some $f\in F$ the intersection $L\cap{\downarrow}f$ has cardinality $\kappa=|L|$, which contradicts the generalized $\K$-Lusin property of $L$.
\smallskip

(2) By the definition of the cardinal $\mathfrak d$, there exists a set $D\subseteq\w^\w$ of cardinality $\mathfrak d$ such that $\w^\w=\bigcup_{f\in D}{\downarrow}f$. Since $L=\bigcup_{f\in D}L\cap{\downarrow}f$ is generalized $\K$-Lusin, for every $f\in D$ the set $L\cap{\downarrow}f$ has cardinality $<\kappa$ and hence $\mathfrak d=|D|\ge\cf(\kappa)$.
\end{proof}

A subset $A$ of a topological space $X$ is said to be {\em concentrated} at a subset $B\subseteq X$ if for every open set $U\subseteq X$ that contains $B$ the complement $A\setminus U$ is countable. Observe that the real line is concentrated at any Bernstein subset of $\IR$.

A subset $A$ of a real line is called {\em $\IQ$-concentrated} if it is concentrated at the set $\IQ$ of rational numbers.

The existence of uncountable $\IQ$-concentrated sets or uncountable $\K$-Lusin sets is independent of \zfc-axioms:

\begin{lemma}\label{l:bQK} The following conditions are equivalent:
\begin{enumerate}
\item[\textup{(1)}] There exists an uncountable $\IQ$-concentrated set.
\item[\textup{(2)}] There exists an uncountable $\K$-Lusin.
\item[\textup{(3)}] $\mathfrak b=\w_1$.
\end{enumerate}
\end{lemma}

\begin{proof}  The equivalence $(1)\Leftrightarrow(2)$ follows from the well-known fact that $\w^\w$ is homeomorphic to $\IR\setminus\IQ$ and the observation that a subset $X\subseteq\IR\setminus\IQ$ is $\K$-Lusin in $\IR\setminus\IQ$ if and only if $X$ is $\IQ$-concentrated.
\smallskip

$(2)\Ra(3)$ Assume that $L$ is an uncountable $\K$-Lusin set. Since subsets of $\K$-Lusin sets are $\K$-Lusin, we can assume that $|L|=\w_1$. To see that $\mathfrak b=\w_1$, it suffices to show that the set $L$ is $\le^*$-unbounded in $\w^\w$. Assuming the converse and applying Lemma~\ref{lem:RemoveTheStar}, we can find a countable set $C\subseteq\w^\w$ such that $L\subseteq\bigcup_{c\in C}{\downarrow}c$. By the Pigeonhole Principle, for some $c\in C$ the set $L\cap {\downarrow}c$ is uncountable, which contradicts the $\K$-Lusin property of $L$ (as the set ${\downarrow}c$ is compact in $\w^\w$).
\smallskip

The implication $(3)\Ra(2)$ follows from Lemma~\ref{l:genL} and the fact that generalized $\K$-Lusin sets of cardinality $\w_1$ are $\K$-Lusin.
\end{proof}

\section{Scattered modifications of topological spaces}

All examples of uncountable FLS spaces constructed in this paper are scattered modifications of separable metrizable spaces. Those modifications are defined as follows.

Let $(X,\tau)$ be a nonempty topological space. A surjective function $h: X\to h[X]$ to an ordinal $h[X]$ is called a {\em height function} if for every $\alpha\in h[X]$ the preimage $h^{-1}\{\alpha\}$ is dense in the set $h^{-1}[{\uparrow}\alpha]$ where ${\uparrow}\alpha\defeq\{\beta\in h[X]:\alpha\le\beta\}$. Every height function $h:X\to h[X]$ determines a topology $\tau_h$ on $X$ generated by the base consisting of the sets
$$B_h(x;U)\defeq\{x\}\cup\{y\in U:h(y)<h(x)\},\quad\mbox{where}\quad x\in U\in\tau.$$
The topological space $X_h=(X,\tau_h)$ is called the {\em scattered $h$-modification} of $X$. Some basic properties of scattered modifications are established in the following propositions.

\begin{proposition}\label{p:Xh1} Let $h:X\to h[X]$ be a height function on a topological space $X$. Then
\begin{enumerate}
\item[\textup{(1)}] $X_h$ is a scattered space.
\item[\textup{(2)}] $X_h^{(\alpha)}=h^{-1}[{\uparrow}\alpha]$ for every ordinal $\alpha\in h[X]$.
\item[\textup{(3)}] The function $h$ coincides with the height function $\hbar$ of $X_h$ and the scattered height $\hbar[X_h]$ of  $X_h$ coincides with the ordinal $h[X]$.
\end{enumerate}
\end{proposition}

\begin{proof} (1) To see that the space $X_h$ is scattered, observe that for every set $A\subseteq X$, every point $a\in A$ with $h(a)=\min h[A]$ is isolated in $A$.
\smallskip

(2) By transfinite induction we shall prove that $X_h^{(\alpha)}=h^{-1}[{\uparrow}\alpha]$ for every ordinal $\alpha\in h[X]$. This equality is trivial for $\alpha=0$. Assume that for some nonzero ordinal $\alpha\in h[X]$ we have proved that $X_h^{(\beta)}=h^{-1}[{\uparrow}\beta]$ for all $\beta\in\alpha$. If the ordinal $\alpha$ is limit, then $$\textstyle X_h^{(\alpha)}=\bigcap_{\beta<\alpha}X_h^{(\beta)}=\bigcap_{\beta<\alpha}h^{-1}[{\uparrow}\beta]=h^{-1}[\bigcap_{\beta<\alpha}{\uparrow}\beta]=h^{-1}[{\uparrow}\alpha].$$
Next, assume that $\alpha=\beta+1$ for some ordinal $\beta$. By the inductive assumption, $X_h^{(\beta)}=h^{-1}[{\uparrow}\beta]$. By the definition of the topology of the space $X_h$, every point of the set $h^{-1}\{\beta\}$ is isolated in $h^{-1}[{\uparrow}\beta]=X_h^{(\beta)}$. Since $h^{-1}\{\beta\}$ is dense in $h^{-1}[{\uparrow}\beta]$ endowed with the subspace topology inherited from $X$, every point of the set $h^{-1}[{\uparrow}\alpha]$ is not isolated in $h^{-1}[{\uparrow}\beta]\subseteq X$ and the definition of the topology $\tau_h$ ensures that it remains non-isolated in $X_h^{(\beta)}=h^{-1}[{\uparrow}\beta]\subseteq X_h$. Therefore, $X^{(\alpha)}_h=(X^{(\beta)}_h)^{(1)}=h^{-1}[{\uparrow}\alpha]$.
\smallskip

(3) The preceding two statements imply that $h=\hbar$ and hence $\hbar[X]=h[X]$.
\end{proof}

Next, we show that the scattered modification inherits some properties of the space.

\begin{proposition}\label{p:Xh2} Let $h:X\to h[X]$ be a height function on a (not necessarily Hausdorff) topological space $X$. 
\begin{enumerate}
\item[\textup{(1)}] If $X$ is first-countable, then so is the space $X_h$.
\item[\textup{(2)}] If $X$ is quasi-metrizable, then so is the space $X_h$.
\item[\textup{(3)}] If $X$ is submetrizable, then so is the space $X_h$.
\item[\textup{(4)}] If $X$ is a regular $T_1$ space, then $X_h$ is almost regular.
\item[\textup{(5)}] If $X$ is hereditarily separable, then $X_h$ is hereditarily separable if and only if for every $\alpha\in h[X]$ the set $h^{-1}\{\alpha\}$ is countable.
\item[\textup{(6)}] If $X$ is first-countable, $h[X]=2$, and $h^{-1}\{1\}$ is countable, then the space $X_h$ has a point-countable base.
\item[\textup{(7)}] If $X$ is a $T_1$ space and the set $h^{-1}\{0\}$ has empty interior in $X$, then the space $X_h$ is not semiregular.
\item[\textup{(8)}] If $X$ is a semiregular $T_1$ space, then the semiregularization of the space $X_h$ coincides with the Bing--Hanner modification $X_M$ of $X$ by the set $M=h^{-1}[{\uparrow}1]$.
\item[\textup{(9)}] If $X$ is metrizable, then  $X_h$ is nearly metrizable if and only if $h^{-1}\{0\}$ is an $F_\sigma$-set in $X$.
\end{enumerate}
\end{proposition}

\begin{proof} (1) If $\mathcal U_x$ is a countable neighborhood base at some point $x\in X$, then $\{B_h(x;U):U\in\U_x\}$ is a countable neighborhood base at $x$ in the space $X_h$.
\smallskip

(2) If $d:X\times X\to \IR_+$ is a quasi-metric generating the topology of $X$, then the function $\rho:X\times X\to\IR_+$ defined by
$$\rho(x,y)=\begin{cases}\min\{1,d(x,y)\}&\mbox{if $y=x$ or $h(y)<h(x)$};\\
1&\mbox{otherwise};
\end{cases}
$$is a quasi-metric generating the topology of the space $X_h$.
\smallskip

(3) Since the identity map $X_h\to X$ is continuous, the submetrizability of $X$ implies the submetrizability of $X_h$.
\smallskip

(4) Assume that $X$ is a regular $T_1$ space. To show that $X_h$ is almost regular, take any regular open set $U$ in $X_h$ and $x\in U$. If $h(x)=0$, then the point $x$ is isolated in  $X_h$ and then the singleton $\{x\}\subseteq U$ is a closed neighborhood of $x$ in $X_h$. So, we assume that $h(x)>0$. By the definition of the topology of the space $X_h$, there exists an open neighborhood $V$ of $x$ in $X$ such that $B_h(x;V)\subseteq U$. Since $X$ is  regular, there exists an open neighborhood $W$ of $x$ in $X$ such that $\overline W\subseteq V$. We claim that the closure of the neighborhood $B_h(x;W)$ in $X_h$ is contained in $U$. Let $z$ be any point in the closure of $B_h(x;W)$ in $X_h$. If $z\in B_h(x;W)$, then $$z\in B_h(x;W)\subseteq B_h(x;V)\subseteq U$$and we are done. So, assume that $z\notin B_h(x;W)$. Since every point of the set $h^{-1}\{0\}$ is isolated in $X_h$, the point $z\in\overline{B_h(x;W)}\setminus B_h(x;W)$ has $h(z)>0$. Since the identity map $X_h\to X$ is continuous, the point $z$ belongs to the closure of the set $B_h(x;W)\subseteq W$ in $X$ and hence $z\in\overline W\subseteq V$. We claim that the set $B_h(z;V)$ is contained in the closure of the set $U$ in $X_h$. Given any point $y\in B_h(z;V)\subseteq V$ and any open neighborhood $O_y$ of $y$ in $X_h$, we should prove that $O_y\cap U\ne \emptyset$. If $h(y)=0$, then $y\in O_y\cap B_h(x;V)\subseteq O_y\cap U$ and we are done. So, assume that $h(y)>0$. By the definition of the topology of the space $X_h$, there exists an open neighborhood $V_y$ of $y$ in $X$ such that $B_h(y;V_y)\subseteq O_y$. By the density of the set $h^{-1}\{0\}$ in $X$, there exists a point $v\in h^{-1}\{0\}\cap V_y\cap V$. Then $v\in B_h(y;V_y)\cap B_h(x;V)\subseteq O_y\cap U$, witnessing that $y\in\overline U$. Therefore, $B_h(z;V)\subseteq\overline U$ and $z$ is an interior point of the closure of $U$ in $X_h$. Since the set $U$ is regular open in $X_h$, $z\in U$ and finally $\overline{B_h(x;W)}\subseteq U$ in $X_h$, witnessing that the space $X_h$ is almost regular.
\smallskip

(5) If the space $X_h$ is hereditarily separable, then for every $\alpha\in h[X]$ the preimage $h^{-1}\{\alpha\}=X_h^{(\alpha)}\setminus X_h^{(\alpha+1)}$ is countable, being a discrete subspace of the hereditarily separable space $X_h$. Now assume that the space $X$ is hereditarily separable and the height function $h$ has countable fibers. Assuming that the space $X_h$ is not hereditarily separable, we can find a nonseparable subspace $Y\subseteq X_h$. Let $h_Y\defeq h{\restriction}_Y$. Consider the sequence of ordinals $(z_\alpha)_{\alpha\in\w_1}$ defined by the recursive formula: $z_0\defeq\min h[Y]$ and $$z_\alpha\defeq\min h[Y\setminus \overline{Y_{<\alpha}}]\quad\mbox{where}\quad \textstyle Y_{<\alpha}\defeq\bigcup_{\beta<\alpha}h_Y^{-1}\{z_\beta\}.$$
Since $Y$ is not separable, for every $\alpha\in\w_1$ the set $Y\setminus\overline {Y_{<\alpha}}$ is not empty and hence the sequence $(z_\alpha)_{\alpha\in\w_1}$ is a well-defined and strictly increasing. For every $\alpha\in\w_1$ choose an element $y_\alpha\in Y\setminus\overline{Y_{<\alpha}}$ such that $h(y_\alpha)=z_\alpha$. Since the space $X$ is hereditarily separable, the subspace $\{y_\alpha:\alpha\in\w_1\}$ of $X$ contains a countable dense subset. Consequently, there exists a  countable ordinal $\beta$ such that the set $D\defeq\{y_\alpha:\alpha<\beta\}$ is dense in $\{y_\alpha:\alpha\in\w_1\}$. Since $y_\beta\notin\overline{Y_{<\beta}}$, there exists a neighborhood $U$ of $y_\beta$ in $X$ such that $B_h(y_\beta;U)\cap Y_{<\beta}=\emptyset$.
By the density of the set $D$ in $\{y_\alpha:\alpha\in\w_1\}\subseteq X$, there exists an ordinal $\alpha<\beta$ such that $y_\alpha\in U$. Observe that $y_\alpha\in Y_{<\beta}$ and $h(y_\alpha)=z_\alpha<z_\beta=h(y_\beta)$. Then $y_\alpha\in Y_{<\beta}\cap B_h(y_\beta;U)$, which contradicts the choice of the neighborhood $U$. This contradiction shows that the space $X_h$ is hereditarily separable.
\smallskip

(6) Assume the space $X$ is first-countable, $h[X]=2$ and the set $h^{-1}\{1\}$ is countable. For every $x\in h^{-1}\{1\}$ fix a countable neighborhood base $\mathcal B_x$ at $x$, and observe that the family
$$\big\{B_h(x;U):x\in h^{-1}\{1\},\;U\in\mathcal B_x\}\cup\{\{x\}:x\in h^{-1}\{0\}\big\}$$is a point-countable base for the space $X_h$.
\smallskip

(7) Assume that $X$ is a $T_1$-space and the set $h^{-1}\{0\}$ has empty interior in $X$. Then $h[X]\ge 2$ and there exists a point $x\in h^{-1}\{1\}$. Assuming that the space $X$ is semiregular, we can find a regular open neighborhood $U\subseteq B_h(x;X)$ of $x$ in $X_h$. By the definition of the topology of the space $X_h$, there exists a neighborhood $V$ of $x$ in $X$ such that $B_h(x;V)\subseteq U$. Since the set $h^{-1}\{0\}$ is dense in the $T_1$ space $X$, the intersection $V\cap h^{-1}\{0\}$ is infinite and the open set $V\setminus\{x\}$ is not empty. Since the set $h^{-1}\{0\}$ has empty interior in $X$, there exists a point $y\in (V\setminus\{x\})\setminus h^{-1}\{0\}$. We claim that $B_{h}(y;V)\subseteq \overline U$ in $X_h$. Take any point $z\in B_{h}(y;V)$. If $h(z)=0$, then $z\in V\cap h^{-1}\{0\}\subseteq B_h(x;V)\subseteq U\subseteq \overline U$ and we are done. If $h(z)>0$, then for every neighborhood $W$ of $z$ in $X$, by the density of the set $h^{-1}\{0\}$ in $X$, we can find a point $w\in W\cap V\cap h^{-1}\{0\}$ and conclude that
$$w\in B_h(z;W)\cap B_h(x;V)\subseteq B_h(z;W)\cap U,$$
which implies that $B_h(z;W)\cap U\ne \emptyset$ and $z\in \overline U$. Then $B_{h}(y;V)\subseteq \overline U$ and hence $y\in U$ as $U$ is regular open. Then $y\in U\subseteq B_{h}(x;X)$ implies $h(y)=0$, which contradicts the choice of $y$. This contradiction shows that the space $X$ is not semiregular.
\smallskip

(8) Assume that the space $X$ is semiregular and $T_1$. Let $\tau_{hsr}$ be the topology of semiregularization of the space $X_h$. We need to show that $\tau_{hsr}$ coincides with the topology $\tau_M$ of the Bing--Hanner modification $X_M$ of $X$ by the set $M=h^{-1}[{\uparrow}1]$. By the semiregularity of $X$, the family $\tau_{ro}$ of  regular open sets in $X$ is a base of the topology of $X$ and the family $\tau_{ro}\cup\{\{x\}:X\setminus M\}$ is a base of the topology of the Bing--Hanner modification $X_M$. The inclusion $\tau_M\subseteq \tau_{hsr}$ will be proved as soon as we check that every set $U\in\tau_{ro}\cup\{\{x\}:x\in X\setminus M\}$ is regular open in $X_h$. If $U=\{x\}$ for some $x\in X\setminus M=h^{-1}\{0\}$, then $\{x\}$ is regular open in $X_h$ since $x$ is an isolated point in $X_h$ and the singleton $\{x\}$ is closed in $X$ and closed-and-open in $X_h$.
If $U$ is regular open in $X$, then  the continuity of the identity map $X_h\to X$ implies that $U$ is open in $X_h$ and its closure $\overline U$ in $X_h$ is contained in the closure $\cl_X(U)$ of $U$ in $X$. It remains to show that $U$ coincides with the interior of $\overline U$ in $X_h$. In the opposite case we can find a point $x\in\overline{U}\setminus U$ and a neighborhood $V$ of $x$ in $X$ such that $B_{h}(x;V)\subseteq\overline{U}$. Observe that the point $x\in\overline{U}\setminus U$ is not isolated in $X_h$ and hence $h(x)>0$. Then the set $V\cap h^{-1}\{0\}$ is dense in $V\subseteq X$ and the inclusion $V\cap h^{-1}\{0\}\subseteq B_{h}(x;V)\subseteq\overline U\subseteq \cl_X(U)$ implies that $V\subseteq \cl_X(V)=\cl_X(V\cap h^{-1}\{0\})\subseteq \cl_X(U)$ and hence $x$ is an interior point of  $\cl_X(U)$ and finally $x\in U$ (because $U$ is regular open in $X$), which contradicts the choice of $x$. This contradiction shows that the set $U$ is regular open in $X_h$ and hence $\tau_M\subseteq\tau_{hsr}$.

To check that $\tau_{hsr}\subseteq\tau_M$, take any regular open set $W$ in $X_h$. To show that $W\in\tau_M$, it suffices to check that every point $x\in W$ is an interior point of the set $W$ in $X_M$. This is clear if $x\notin M$. So, assume that $x\in M$ and hence $h(x)>0$. By the definition of the topology of the space $X_h$, there exists a neighborhood $U_x$ of $x$ in $X$ such that $B_h(x;U_x)\subseteq W$. We claim that $U_x$ is contained in the closure of the set $W$ in $X_h$. Given any point $y\in U_x$ and neighborhood $V_y$ of $y$ in $X_h$, we need to show that $V_y\cap W\ne\emptyset$. If $h(y)<h(x)$, then $y\in V_y\cap B_h(x;U_x)\subseteq V_y\cap W$ and we are done. So, assume that $h(y)\ge h(x)>0$.
By the definition of the topology of the space $X_h$, there exists a neighborhood $U_y$ of $y$ in $X$ such that $B_h(y;U_y)\subseteq V_y$. By the density of the set $h^{-1}\{0\}$ in $X$, the intersection $h^{-1}\{0\}\cap U_y\cap U_x$ contains some point $z$. Then $z\in B_h(y;U_y)\cap B_h(x;U_x)\subseteq V_y\cap W$ and hence $V_y\cap W\ne\emptyset$ and $U_x\subseteq\overline{W}$. The continuity of the identity map $X_h\to X$ ensures that the set $U_x$ remains open in the space $X_h$ and hence $U_x$ is contained in the interior of the set $\overline{W}$ which is equal to $W$ (as $W$ is regular open in $X_h$). Therefore, $x\in U_x\subseteq W$ is an interior point of the set $W$ in $X_M$ and $W\in\tau_M$. 
\smallskip

(9) Assume that the space $X$ is metrizable. By the preceding statement, the semiregularization $X_{hsr}$ of the space $X_h$ coincides with the Bing--Hanner  modification $X_M$ of $X$ by its subset $M=h^{-1}[{\uparrow}1]$. By Lemma~\ref{l:Michael}, the space $X_M=X_{hsr}$ is metrizable if and only if the set $h^{-1}\{0\}=X\setminus M$ is of type $F_\sigma$ in $X$.
\end{proof}

Finally we find some conditions guaranteeing that the scattered modification of a topological space is Lindel\"of.

We recall that a topological space $X$ is {\em concentrated} at a subset $B\subseteq X$  if every open set containing $B$ has countable complement in $X$.

\begin{lemma}\label{l:XhL1} Let $X$ be a topological space and $h:X\to h[X]$ be a height function. The scattered modification $X_h$ of $X$ is Lindel\"of if there exists an ordinal $\alpha\in h[X]$ such that the subspace $h^{-1}[{\uparrow}\alpha]$ of $X_h$ is Lindel\"of and $X$ is concentrated at $h^{-1}[{\uparrow}\alpha]$.
\end{lemma}

\begin{proof} Let $\U$ be any open cover of the space $X_h$. For every $x\in X$ find an open neighborhood $U_x\in\U$ of $x$ and an open neighborhood $V_x$ of $x$ in $X$ such that $B_h(x;V_x)\subseteq U_x$. Since the subspace $h^{-1}[{\uparrow}\alpha]$ of $X_h$ is Lindel\"of, there exists a countable subset $C\subseteq h^{-1}[{\uparrow}\alpha]$ such that $h^{-1}[{\uparrow}\alpha]\subseteq\bigcup_{x\in C}B_h(x;V_x)$. Consider the open set $V\defeq\bigcup_{x\in C}V_x$ in $X$ and observe that  
$$h^{-1}[{\uparrow}\alpha]\subseteq\bigcup_{x\in C}B_h(x;V_x)=\bigcup_{x\in C}V_x=V.$$
The equality $\bigcup_{x\in C}B_h(x;V_x)=\bigcup_{x\in C}V_x$ follows from the observation that $V_x\setminus h^{-1}[{\uparrow}\alpha]=B_h(x;V_x)\setminus h^{-1}[{\uparrow}\alpha]$ for every $x\in C\subseteq h^{-1}[{\uparrow}\alpha]$.

Since $X$ is concentrated at $h^{-1}[{\uparrow}\alpha]$, the complement $X\setminus V$ is countable. It is easy to see that $\{U_x:x\in C\cup(X\setminus V)\}\subseteq\U$ is a countable cover of $X$, witnessing that the space $X_h$ is Lindel\"of.
\end{proof}

\begin{lemma}\label{l:XhL2} Let $X$ be a topological space and $h:X\to h[X]$ be a height function. The scattered modification $X_h$ of $X$ is Lindel\"of if there exists a strictly increasing sequence of ordinals $(\alpha_n)_{n\in\w}$ in $h[X]$ satisfying the following conditions:
\begin{enumerate}
\item[\textup{(1)}] for the ordinal $\alpha=\sup_n\alpha_n$, the subspace $h^{-1}[{\uparrow}\alpha]$ of $X_h$ is Lindel\"of;
\item[\textup{(2)}] $\alpha_0=0$ and for every $n\in\w$ the subspace $h^{-1}[{\uparrow}\alpha_n\cap{\downarrow}\alpha_{n+1}]$ of $X$ is concentrated at $h^{-1}\{\alpha_{n+1}\}$ and the subspace $h^{-1}\{\alpha_{n+1}\}$ of $X$ is Lindel\"of.
\end{enumerate}
\end{lemma}

\begin{proof} Let $\U$ be any open cover of the space $X_h$. For every $x\in X$ find an open neighborhood $U_x\in\U$ of $x$ and an open neighborhood $V_x$ of $x$ in $X$ such that $B_h(x;V_x)\subseteq U_x$. Since the subspace $h^{-1}[{\uparrow}\alpha]$ of $X_h$ is Lindel\"of, there exists a countable subset $C\subseteq h^{-1}[{\uparrow}\alpha]$ such that $h^{-1}[{\uparrow}\alpha]\subseteq\bigcup_{x\in C}U_x$. By our assumption, for every $n\in\w$ the subspace $h^{-1}\{\alpha_{n+1}\}$ of $X$ is Lindel\"of. Consequently, there exists a countable set $C_n\subseteq h^{-1}\{\alpha_{n+1}\}$  such that $h^{-1}\{\alpha_{n+1}\}\subseteq\bigcup_{x\in C_n}V_x$. Since $h^{-1}[{\uparrow}\alpha_n\cap{\downarrow}\alpha_{n+1}]$ is concentrated at $h^{-1}\{\alpha_{n+1}\}$, the set $$D_n=(h^{-1}[{\uparrow}\alpha_n\cap{\downarrow}\alpha_{n+1}])\setminus\bigcup_{x\in C_n}V_x$$is countable.
We claim that $\U'=\{U_x:x\in C\cup\bigcup_{n\in\w}(C_n\cup D_n)\}\subseteq\U$ is a countable cover of $X$. Take any point $y\in X$. If $y\in h^{-1}[{\uparrow}\alpha]$, then $y\in\bigcup_{x\in C}U_x$ by the choice of $C$. If $y\notin h^{-1}[{\uparrow}\alpha]$, then there exists a unique number $n\in\w$ such that $\alpha_n\le h(y)<\alpha_{n+1}$. If $y\in D_n$, then $y\in U_y\in\U'$. If $y\notin D_n$, then $y\in V_x$ for some $x\in C_n\subseteq h^{-1}\{\alpha_{n+1}\}$ and then $y\in B_h(x;V_x)\subseteq U_x\in \U'$.

Therefore, $\U'$ is a countable subcover of $\U$, witnessing that the space $X_h$ is Lindel\"of.
\end{proof}

\begin{definition} A scattered topological space will be called {\em special} if it is homeomorphic to a scattered modification of some metrizable  space.
\end{definition}

Proposition~\ref{p:Xh2} implies the following fact.

\begin{corollary} If a scattered space $X$ is special (and separable), then $X$ is almost regular, first-countable, submetrizable, quasi-metrizable (and nearly metrizable).
\end{corollary}

\section{Constructions of special FLS spaces}\label{s:constructions}

In this section we apply special FLS spaces to construct FLS spaces of given cardinality and height. In the following four theorems we present four such constructions.

\begin{theorem}\label{thm:Existence} For any ordinal $\lambda<\mathfrak c^+$ with $\cf(\lim\lambda)=\w$, there exists a separable special FLS space of  cardinality $\continuum$ and scattered height $\lambda$.
\end{theorem}

\begin{proof} By Lemma~\ref{lem:Bernstein}, there is a function $h:\IR\to\lambda$ such that  
\begin{itemize}
\item for every ordinal $\alpha\in\{0\}\cup [\lim\lambda,\lambda)$ the preimage $h^{-1}\{\alpha\}$ is countable and dense in $\mathbb R$;
\item  for every nonzero ordinal $\alpha\in\lim \lambda$, the preimage $h^{-1}\{\alpha\}$ is a Bernstein set in $\IR$.
\end{itemize}

Since each fiber of the function $h$ is dense in $\IR$, $h$ is a height function. By Proposition~\ref{p:Xh2}, the scattered modification $\IR_h$ of the real line is a scattered first-countable space of scattered height $\lambda$. Since the set $h^{-1}\{0\}$ of isolated points of $\IR_h$ is countable, the space $\IR_h$ is separable.

 By  Lemma~\ref{l:XhL2}, the space $\IR_h$ is Lindel\"of. So, $\IR_h$ is a separable special FLS space of cardinality $\mathfrak c$ and scattered height $\lambda$.
\end{proof}

\begin{theorem}\label{t:construct2} Under $\cov{M}=\mathfrak d$, for every successor ordinal $\lambda<\mathfrak c^+$ with $\cf(\lim \lambda)=\cf(\mathfrak d)$ there exists a separable special FLS space $X$ of cardinality $\mathfrak c$ and scattered height $\lambda$.
\end{theorem}

\begin{proof} Let $\kappa\defeq\cf(\lim\lambda)=\cf(\mathfrak d)$. Choose a strictly increasing sequence of ordinals $(\lambda_\alpha)_{\alpha\in\kappa}$ such that 
\begin{itemize}
\item $\lambda_0=0$; 
\item $\lambda_\alpha=\sup_{\beta\in\alpha}\lambda_\beta$ for every infinite limit ordinal $\alpha\in\kappa$;
\item $\sup_{\alpha\in\kappa}\lambda_\alpha=\lim\lambda$;
\item $\cf(\lambda_{\alpha+1})=\w$ for every ordinal $\alpha\in\kappa$.
\end{itemize}

Using Lemma~\ref{lem:Bernstein}, construct a family of sets $(B_\beta)_{\beta\in \lim\lambda}$ such that 
\begin{itemize}
\item the set $B_0$ is countable and dense in $\w^\w$;
\item for every ordinal $\alpha\in[1,\lim\lambda)$ the set $B_\alpha$ is Bernstein in $\w^\w$;
\item for every $\alpha\in\kappa$ the family $(B_\beta)_{\beta\in [\lambda_\alpha,\lambda_{\alpha+1})}$ is a disjoint cover of $\w^\w$. 
\end{itemize}

\begin{claim}\label{cl:S} There exists a sequence $(S_\alpha)_{\alpha\in\mathfrak d}$ of $\sigma$-compact sets in $\w^\w$ such that
\begin{enumerate}
\item[\textup{(1)}] $B_0\subseteq S_0$;
\item[\textup{(2)}] every compact subset of $\w^\w$ is contained in some set $S_\alpha$;
\item[\textup{(3)}] for every $\alpha\in\mathfrak d$ and every Bernstein set $B$ in $\w^\w$ the intersection $B\cap S_\alpha\setminus \bigcup_{\beta\in\alpha}S_\beta$ is dense in $\w^\w$ and has cardinality of continuum.
\end{enumerate}
\end{claim}

\begin{proof} By definition of $\mathfrak d$, there exists a transfinite sequence $\{g_\alpha\}_{\alpha\in\mathfrak d}\subseteq\w^\w$ such that $\w^\w=\bigcup_{\alpha\in\kappa}{\downarrow}g_\alpha$. Since every countable set in $\w^\w$ is $\le^*$-bounded, we can replace $g_0$ by a larger function and assume that $x\le^* g_0$ for every $x\in B_0$.

Using the definition of the cardinal $\mathfrak d$, construct inductively a transfinite sequence $(f_\alpha)_{\alpha\in\mathfrak d}$ in $\IN^\w$ such that 
for every $\alpha\in\mathfrak d$ the following conditions are satisfied:
\begin{itemize}
\item[(i)] $g_\alpha+1\le f_\alpha$;
\item[(ii)] $f_\alpha\not\le^* f_\beta+2$ for every $\beta\in\alpha$.
\end{itemize}
For every $\alpha\in\mathfrak d$, consider the subsets
$$S_\alpha\defeq\{f\in\w^\w:f\le^* f_\alpha\}\quad\mbox{and}\quad S_{<\alpha} \defeq \textstyle \bigcup_{\beta\in\alpha}S_\beta$$
of $\w^\w$. 
The choice of $g_0\le f_0$ ensures that $B_0\subseteq S_0$. The choice of the transfinite sequence $(g_\alpha)_{\alpha\in\mathfrak d}$ guarantees that  for every compact subset $K$ of $\w^\w$ there exists $\alpha\in\mathfrak d$ such that $K\subseteq{\downarrow}g_\alpha\subseteq {\downarrow}f_\alpha\subseteq S_\alpha$. To show that for every $\alpha\in\mathfrak d$ and every  Bernstein set $B$ in $\w^\w$ the intersection $B\cap S_\alpha\setminus S_{<\alpha}$ is dense in $\w^\w$ and has cardinality of continuum, it suffices to check that for every nonempty open set $U\subseteq\w^\w$ the intersection $U\cap S_\alpha\setminus S_{<\alpha}$ contains an uncountable compact subset.

Given a nonempty open set $U\subseteq\w^\w$, fix any function $u\in U$ and find $m\in\w$ such that any function $v\in\w^\w$ with $v{\restriction}_m=u{\restriction}_m$ belongs to $U$. Let $V$ be the set of all functions $v\in\w^\w$ that satisfy the following conditions:
\begin{itemize}
\item $v(n)=u(n)$ for every $n<m$;
\item $v(n)\in \{f_\alpha(n), f_\alpha(n)-1\}$ for every $n\ge m$.
\end{itemize}
It is clear that the set $V$ is compact and uncountable. It is also clear that $V\subseteq U\cap {\downarrow}f_\alpha\subseteq S_\alpha$. The condition (ii) of the inductive construction ensures that $V$ is disjoint with the set $S_{<\alpha}$.
\end{proof}

Let $(d_\alpha)_{\alpha\in\kappa}$ be a strictly increasing sequence of ordinals such that $d_0=0$ and $\sup_{\alpha\in\kappa}d_\alpha=\mathfrak d$.

For any ordinal $\beta\in\lim\lambda$, find a unique ordinal $\alpha\in\kappa$ such that $\beta\in[\lambda_\alpha,\lambda_{\alpha+1})$ and put 
$$D_\beta\defeq B_\beta\cap 
S_{<d_{\alpha+1}}\setminus S_{<d_{\alpha}}$$
where, as in the proof of the previous claim, $S_{<\beta}\defeq\bigcup_{\gamma\in\beta}S_{\gamma}$ for an ordinal $\beta\in\mathfrak d$. Observe that the set $D_0=B_0\cap S_{<d_1}=B_0$ is countable. 
Claim~\ref{cl:S}(3) ensures that for every nonzero ordinal $\beta\in\lim\lambda$ the set $D_\beta$ is dense in $\w^\w$ and has cardinality $\mathfrak c$. Therefore, $(D_\beta)_{\beta\in\lim\lambda}$ is a disjoint cover of $\w^\w$ by dense sets such that $D_0$ is countable. 

Since the space $\w^\w$ is homeomorphic to the space $\IR\setminus\IQ$ of irrational numbers, there exists a metrizable compactification $X$ of $\w^\w$ such that $Q=X\setminus\w^\w$ is a countable dense set in $X$. Let $(D_\beta)_{\beta\in[\lim\lambda,\lambda)}$ be a finite disjoint cover of $Q$ by countable dense subsets of $X$. 

Let $h:X\to\lambda$ be the unique function such that $h^{-1}\{\beta\}=D_\beta$ for every $\beta\in\lambda$. The density of the sets $D_\beta$ in $X$ implies that $h$ is a height function on $X$ and hence the scattered $h$-modification $X_h$ of the space $X$ is a special first-countable scattered space of cardinality  $|X_h|=|X|=\mathfrak c$ and height $\lambda$, according to Propositions~\ref{p:Xh1} and \ref{p:Xh2}. Since the set $h^{-1}\{0\}=D_0=B_0$ of isolated points of $X_h$ is countable, the scattered space $X_h$ is separable. It remains to prove that the space $X_h$ is Lindel\"of. 

To derive a contradiction, assume that $X_h$ is not Lindel\"of and find an open cover  $\U$ of $X_h$ that has no countable subcovers. For every point $x\in X$ find a set $U_x\in\U$ containing $x$, and an open neighborhood $V_x$ of $x$ in $X$ such that $B_h(x;V_x)\subseteq U_x$. 

\begin{claim}\label{cl:interior} Under $\cov{M} = \mathfrak d$, for every nonempty compact set $K$ in $\w^\w$ there exists a countable set $C\subset  \w^\w$ such that the intersection $K\cap\bigcup_{x\in C}U_x$ has nonempty interior in $K$.
\end{claim}

\begin{proof}  If $K$ has an isolated point $x$, then put $C=\{x\}$ and finish the proof. So, we assume that $K$ has no isolated points. Let $\gamma\in\mathfrak d$ be the smallest ordinal such that the intersection $K\cap S_\gamma$ is nonmeager in $K$. By Claim~\ref{cl:S}(2), the ordinal $\gamma$ is well-defined. Since the set $S_\gamma$ is $\sigma$-compact, we can apply the Baire Theorem and find a nonempty open set $W$ in $K$ such that $\overline{W}\subseteq S_\gamma$.  The minimality of $\gamma$ ensures that for every $\beta\in \gamma$ the intersection $K\cap S_\beta$ is meager in $K$.

Find a unique ordinal $\alpha\in\kappa$ such that $d_\alpha\le\gamma<d_{\alpha+1}$. Since the ordinal $\lambda_{\alpha+1}$ has countable cofinality, there exists a strictly increasing sequence of ordinals $(\lambda_{\alpha,n})_{n\in\w}$ such that $\lambda_\alpha<\lambda_{\alpha,0}$ and $\sup_{n\in\w}\lambda_{\alpha,n}=\lambda_{\alpha+1}$. Since the space $\w^\w$ is hereditarily Lindel\"of, for every $n\in\IN$ there exist a countable set $C_n\subseteq D_{\lambda_{\alpha,n}}$ such that $D_{\lambda_{\alpha,n}}\subseteq \bigcup_{x\in C_n}V_x$. We claim that the set $\Gamma_n\defeq \overline W\setminus (S_{<d_\alpha}\cup\bigcup_{x\in C_n}V_x)$ is countable. To derive a contradiction, assume that the set $\Gamma_n$ is uncountable.
Since every scattered subspace of the compact metrizable space $K$ is countable, the crowded part $\Gamma_n^{(\infty)}$ of $\Gamma_n$ is uncountable. Then the closure $\hat\Gamma_n$ of $\Gamma_n^{(\infty)}$ in $K$ is a crowded compact subspace of $\w^\w$. Fix any countable dense subset $\ddot\Gamma_n$ in $\Gamma_n^{(\infty)}$. By the Aleksandrov--Urysohn Theorem \cite[7.7]{Ke}, the complement $\check\Gamma_n\defeq \hat\Gamma_n\setminus \ddot\Gamma_n$ is homeomorphic to $\w^\w$. Observe that for every $\beta\in d_\alpha$ the intersection $\hat\Gamma_n\cap S_\beta$ is a $\sigma$-compact set, disjoint with $\Gamma_n\supset \ddot\Gamma_n$ and hence $\hat\Gamma_n\cap S_\beta$ is a $\sigma$-compact set in $\check\Gamma_n$. By Lemma~\ref{l:d}, the union $\bigcup_{\beta\in d_\alpha}(\hat\Gamma_n\cap S_{\beta})=\check\Gamma_n\cap S_{<d_\alpha}$ is disjoint with some uncountable compact subset $\Lambda$ of $\check \Gamma_n\subseteq\hat\Gamma_n$. Observe that
$$\Lambda\subseteq\hat\Gamma_n\subseteq \overline{\Gamma_n}\subseteq\overline{W}\setminus\textstyle\bigcup_{x\in C_n}V_x$$
 and hence 
 $$\textstyle \Lambda\subseteq \big(\overline W\setminus\bigcup_{x\in C_n}V_x\big)\setminus S_{<d_\alpha}\subseteq (S_\gamma\setminus S_{<d_\alpha})\setminus\bigcup_{x\in C_n}V_x$$
and
$$\textstyle B_{\lambda_{\alpha,n}}\cap\Lambda\subseteq B_{\lambda_{\alpha,n}}\cap\big(S_\gamma\setminus S_{<d_\alpha}\big)\setminus\bigcup_{x\in C_n}V_x\subseteq D_{\lambda_{\alpha,n}}\setminus\bigcup_{x\in C_n}V_x=\emptyset,
$$which is not possible as $B_{\lambda_{\alpha,n}}$ is a Bernstein set in $\w^\w$ and $\Lambda$ is an uncountable compact set in $\w^\w$. This contradiction shows that the set $\Gamma_n$ is countable. The definition of the set $\Gamma_n$ ensures that 
$$\textstyle \overline{W}\subseteq S_{<d_\alpha}\cup\Gamma_n\cup \bigcup_{x\in C_n}V_x.$$

In particular, $\overline{W}\subseteq S_{<d_\alpha}\cup\Gamma_1\cup\bigcup_{x\in C_1}V_x$. By the choice of the ordinal $\gamma$, for every $\beta\in d_\alpha\le\gamma$ the set $W\cap S_\beta$ is meager in $W$. Since $|d_\alpha|<\mathfrak d=\cov{M}$, the set $W\cap (S_{<d_\alpha}\cup\Gamma_1)$ has empty interior in $W$ and hence the open set $V\defeq W\cap\bigcup_{x\in C_1}V_x$ is dense in $W$. In particular, $V$ is a nonempty open set in $K$.

Consider the countable set $C\defeq\bigcup_{n\in\w}(\Gamma_n\cup C_n)$. We claim that $V\subseteq \bigcup_{x\in C}U_x$. Given any $v\in V$, we should find $x\in C$ such that $v\in U_x$. If $v\in C$, then we can put $x=v$. So, assume that $v\notin C$. 

If $v\in S_{<d_\alpha}$, then $h(v)<\lambda_{\alpha}$. Since $v\in V\subseteq \bigcup_{x\in C_1}V_x$, there exists $x\in C_1\subseteq C$ such that $v\in V_x$. Taking into account that $v\in V_x$ and  $h(v)<\lambda_\alpha\le h(x)$, we conclude that $v\in B_h(x;V_x)\subseteq U_x$. Next, assume that $v\notin S_{<d_\alpha}$. In this case $v\in V\setminus S_{<d_\alpha}\subseteq W\setminus S_{<d_\alpha}\subseteq S_\gamma\setminus S_{<d_\alpha}\subseteq S_{<d_{\alpha+1}}\setminus S_{<d_\alpha}=\bigcup_{\beta\in[\lambda_\alpha,\lambda_{\alpha+1})}D_\beta$ and hence $h(v)\in [\lambda_\alpha,\lambda_{\alpha+1})$. 
Since $\lim_{n\in\w}\lambda_{\alpha,n}=\lambda_{\alpha+1}$, there exists $n\in\w$ such that $h(v)<\lambda_{\alpha,n}$.

Since 
$$\textstyle v\in (\overline{W}\setminus S_{<d_\alpha})\setminus C\subseteq \big(S_{<d_\alpha}\cup\Gamma_n\cup\bigcup_{x\in C_n}V_x\big)\setminus (S_{<d_\alpha}\cup C)\subseteq \bigcup_{x\in C_n}V_x,$$ there exists $x\in C_n\subseteq D_{\lambda_{\alpha,n}}\cap C$ such that $v\in V_x$. It follows from $h(v)<\lambda_{\alpha,n}=h(x)$ that $v\in B_h(x;V_x)\subseteq U_x$.

Therefore, the set $K\cap \bigcap_{x\in C}U_x\supseteq V$ has nonempty interior in $K$, which completes the proof of Claim~\ref{cl:interior}.
\end{proof}

Observe that $K_0\defeq X\setminus\bigcup_{x\in Q}V_x$ is a compact subset of $X\setminus Q=\w^\w$. Since $X=K_0\cup \bigcup_{x\in Q}V_x=K_0\cup\bigcup_{x\in Q}B_h(x;V_x)\subseteq K_0\cup\bigcup_{x\in Q}U_x$ and the cover $\U$ has no countable subcovers, no countable subfamily of $\U$ covers the set $K_0$. We shall inductively construct a strictly increasing sequence $(W_\alpha)_{\alpha\in\w_1}$ of open sets in $K_0$ and a sequence $(C_\alpha)_{\alpha\in\w_1}$ of countable sets in $\w^\w$ such that $W_\alpha\subseteq\bigcup_{x\in C_\alpha}U_x$ for every $\alpha\in\w_1$. We start the inductive construction by setting $W_0=C_0=\emptyset$. Assume that for some nonzero ordinal $\alpha\in\w_1$ we have constructed a strictly increasing sequence $(W_\beta)_{\beta\in\alpha}$ of open sets in $K_0$ and a sequence $(C_\beta)_{\beta\in\alpha}$ of countable sets in $\w^\w$ such that $W_\beta\subseteq\bigcup_{x\in C_\beta}U_x$ for every $\beta\in\alpha$. Consider the compact set $K_\alpha=K_0\setminus \bigcup_{\beta\in\alpha}W_\alpha$ and observe that $K_0=K_\alpha\cup\bigcup_{\beta\in\alpha}W_\beta\subseteq K_\alpha\cup\bigcup_{\beta\in \alpha}\bigcup_{x\in C_\beta}U_x$. Since no countable subfamily of $\U$ covers $K_0$, no countable subfamily of $\U$ covers $K_\alpha$.  This implies that the compact set $K_\alpha$ is not empty. By Claim~\ref{cl:interior}, there exists a countable set $C_\alpha'\subset \w^\w$ such that the set $K_\alpha\cap \bigcup_{x\in C_\alpha'}U_x$ has nonempty interior in $K_\alpha$. Consequently, there exists an open set $W_\alpha'$ in $K_0$ such that 
$$\textstyle \emptyset\ne W'_\alpha\cap K_\alpha\subseteq K_\alpha\cap\bigcup_{x\in C_\alpha'}U_x.$$ Put $$\textstyle W_\alpha\defeq W'_\alpha\cup\bigcup_{\beta\in \alpha}W_\beta=W_\alpha'\cup (K\setminus K_\alpha)\quad\mbox{and}\quad C_\alpha\defeq C_\alpha'\cup\bigcup_{\beta\in\alpha}C_\beta.$$

After completing this inductive construction, we obtain a strictly increasing sequence $(W_\alpha)_{\alpha\in\w_1}$ of open sets in $K_0$, which contradicts the hereditary Lindel\"of property of the compact metrizable space $K_0$. This contradiction shows that the space $X_h$ is Lindel\"of.
\end{proof}

Note that the hypothesis $\cov{M} = \mathfrak d$ was used in the preceding proof only to establish Claim~\ref{cl:interior}. However, it seems plausible that a more careful choice of the sets $S_\a$ could eliminate this hypothesis, and lead to a proof of the preceding theorem from $\mathsf{ZFC}$ alone.

\begin{question}
Is there a sequence $(S_\alpha)_{\alpha < \kappa}$ of compact subsets of $\w^\w$ such that 
\begin{enumerate}
\item[\textup{(1)}] every compact subset of $\w^\w$ is contained in some $S_\a$, and 
\item[\textup{(2)}] for any compact $K \subseteq \w^\w$, if $S_{<\a} \cap K$ is non-meager in $K$ then for some $\b < \a$ the set $S_\b \cap K$ is non-meager in $K$?
\end{enumerate} 
\end{question}
\noindent Such a sequence is sufficient to carry out the above proof. The hypothesis $\cov{M} = \mathfrak{d}$ implies a positive answer to this question: it implies that $({\downarrow}f_\a)_{\a\in \mathfrak d}$ is such a sequence whenever $\set{f_\a}{\a \in \mathfrak d}$ is a $\leq$-dominating family. But we do not know whether a negative answer to the question is also consistent.

Next, we construct uncountable FLS spaces using the existence of uncountable $\K$-Lusin sets in $\w^\w$.

\begin{theorem}\label{t:construct3} If there exists a $\K$-Lusin set of uncountable cardinality $\kappa$, then for every infinite ordinal $\lambda<\mathfrak c^+$ with $\cf(\lambda)\le\w$ and $\cf(\lim\lambda)\le\kappa$ there exists a separable special FLS space of cardinality $\mathfrak c$ and scattered height $\lambda$.
\end{theorem}

\begin{proof} 
Assume that there exists a $\K$-Lusin set of cardinality $\kappa$. 
 Fix an infinite ordinal $\lambda<\mathfrak c^+$ with $\cf(\lambda)\le\w$ and $\cf(\lim\lambda)\le\kappa$.  If $\cf(\lim\lambda)=\w$, then a special FLS space of cardinality $\mathfrak c$ and scattered height $\lambda$ exists by Theorem~\ref{thm:Existence}. So, we assume that $\cf(\lim\lambda)>\w$, which implies that  the ordinal $\lambda$ is successor and hence $\lim\lambda<\lambda$. Since subsets of $\K$-Lusin sets are $\K$-Lusin, we lose no generality assuming that $\kappa=\cf(\lim\lambda)>\w$. Choose a strictly increasing sequence of ordinals $(\lambda_\alpha)_{\alpha\in\kappa}$ having the following properties:
 \begin{itemize}
 \item $\lambda_0=0$;
 \item $\lim\lambda=\sup_{\alpha\in\kappa}\lambda_\alpha$;
 \item $\lambda_\alpha=\sup_{\beta\in\alpha}\lambda_\beta$ for every infinite limit ordinal $\alpha\in\kappa$;
 \item $\cf(\lambda_{\alpha+1})=\w$ for every $\alpha\in\kappa$.
 \end{itemize} 
 Consider the unique function $\xi:\lambda\to\kappa+1$ such that $\xi^{-1}\{\kappa\}=[\lim\lambda,\lambda)$ and $\xi^{-1}\{\alpha\}=[\lambda_\alpha,\lambda_{\alpha+1})$ for every $\alpha\in\kappa$. 
 
For a function $f\in\w^\w$, let ${\Uparrow}f\defeq\{g\in\w^\w:{f\le^* g}\}$. It is clear that ${\Uparrow}f$ is an $F_\sigma$-subset of $\w^\w$.
 
\begin{claim}\label{cl1} There exist a $\K$-Lusin set $\{f_\alpha\}_{\alpha\in\kappa}\subseteq\w^\w$ such that for every Bernstein set $B$ in $\w^\w$ and every $\alpha\in\kappa$ the set
$$\textstyle B\cap\big({\Uparrow}f_\alpha\setminus\bigcup_{\beta\in\alpha}{\Uparrow}f_\beta\big)$$is dense in $\w^\w$ and has cardinality $\mathfrak c$.
\end{claim}

\begin{proof} By our assumption, there exists a $\K$-Lusin set $L\subseteq\w^\w$ of cardinality $\kappa$. Write $L$ as $\{g_\alpha\}_{\alpha\in\kappa}$ where the functions $g_\alpha$ are pairwise distinct. By \cite[2.1]{Hodel}, there exists a family $(A_\alpha)_{\alpha\in\kappa}$ of infinite subsets of $\w$ such that $A_\alpha\cap A_\beta$ is finite for any distinct ordinals $\alpha,\beta\in \kappa$. For every ordinal $\alpha\in\kappa$, let $f_\alpha\in\w^\w$ be the unique function such that
\begin{itemize}
\item  $f_\alpha(2n)=g_\alpha(n)$ for every $n\in\w$;
\item $f_{\alpha}(2n+1)=2$ for every $n\in A_\alpha$;
\item $f_\alpha(2n+1)=0$ for every $n\in\w\setminus A_\alpha$.
\end{itemize}
The $\K$-Lusin property of the set $\{g_\alpha\}_{\alpha\in\kappa}$ implies the $\K$-Lusin property of the set $\{f_\alpha\}_{\alpha\in\kappa}$ in $\w^\w$.

It remains to check that for every $\alpha\in\kappa$ and every Bernstein set $B$ in $\w^\w$ the intersection\newline $B\cap{\Uparrow} f_\alpha\setminus\bigcup_{\beta<\alpha}{\Uparrow} f_\beta$ has cardinality of continuum and is dense in $\w^\w$.

Since $B$ has nonempty intersection with each uncountable compact subset of $\w^\w$ and each uncountable compact set in $\w^\w$ contains continuum many pairwise disjoint uncountable compact subsets, it suffices to show that for every non-empty open set $U\subseteq \w^\w$ the intersection $U\cap {\Uparrow} f_\alpha\setminus\bigcup_{\beta<\alpha}{\Uparrow} f_\beta$ contains an uncountable compact subset. Fix any function $u\in U$ and find $m\in\w$ such that any function $v\in\w^\w$ with $v{\restriction}_m=u{\restriction}_m$ belongs to $U$. Let $V$ be the set of all functions $v\in\w^\w$ that satisfy the following conditions:
\begin{itemize}
\item $v(n)=u(n)$ for any $n<m$;
\item $v(n)= f_\alpha(n)$ for any even $n\ge m$;
\item $v(2n+1)=2$ for any $n\in A_\alpha$ with $2n+1\ge m$;
\item $v(2n+1)\in\{0,1\}$ for any  $n\in\w\setminus A_\alpha$ with $2n+1\ge m$.
\end{itemize}
It is clear that the set $V$ is compact and uncountable. It is also clear that $V\subseteq U\cap {\Uparrow} f_\alpha$. Assuming that $V$ has a common point with the set $\bigcup_{\beta<\alpha}{\Uparrow} f_\beta$, we can find an ordinal $\beta<\alpha$ and a function $v\in V\cap{\Uparrow} f_\beta$. 
Since $f_\beta\le^* v$, there exists $l\ge m$ such that $f_\beta(n)\le v(n)$ for all $n\ge l$. 
Since the set $A_\beta\setminus A_\alpha$ is infinite, there exists $n\in A_\beta\setminus A_\alpha$ such that $2n+1\ge l\ge m$. Then $v(2n+1)\le 1<2=f_\beta(2n+1)\le v(2n+1)$, which is a desired contradiction showing that $V\subseteq U\cap{\Uparrow} f_\alpha\setminus\bigcup_{\beta<\alpha}{\Uparrow} f_\beta$.
\end{proof} 

Let $\{f_\alpha\}_{\alpha\in\kappa}$ be the $\K$-Lusin set from Claim~\ref{cl1}. 
 By Lemma~\ref{lem:Bernstein}, there exists a family $(B_\alpha)_{\alpha\in\lim\lambda}$ of pairwise disjoint Bernstein sets in $\w^\w$. Since $\w^\w$ is homeomorphic to the space of irrational numbers, there exists a metrizable compactification $K$ of $\w^\w$ whose remainder $Q=K\setminus\w^\w$ is countable. By Claim~\ref{cl1}, the set $B_0\cap {\Uparrow}f_0$ is dense in $\w^\w$ and hence contains a countable dense set $D_0$. By Claim~\ref{cl1}, for every nonzero ordinal $\alpha\in\lim\lambda$ the set 
 $$\textstyle D_\alpha\defeq B_\alpha\cap{\big(}{\Uparrow} f_{\xi(\alpha)}\setminus\bigcup_{\gamma<\xi(\alpha)}{\Uparrow} f_\gamma{\big)}.$$
has cardinality of continuum and is dense in $\w^\w$.
 
Let $\{D_\alpha\}_{\alpha\in[\lim\lambda,\lambda)}$ be a (finite) cover of $Q=K\setminus \w^\w$ by pairwise disjoint dense sets in $Q$. It follows that $(D_\alpha)_{\alpha\in\lambda}$ is a family of pairwise disjoint dense sets in the compact space $K$.
 
Let $X\defeq\bigcup_{\alpha\in\lambda}D_\alpha\subseteq K$ and $h:X\to\lambda$ be the height function such that $h^{-1}\{\alpha\}=D_\alpha$ for every $\alpha\in\lambda$.
By Proposition~\ref{p:Xh1}, the space $X_h$ is first-countable and scattered with scattered height $\lambda$ and cardinality $|X|=|{D_1}|=\mathfrak c$. Since the set $D_0$ of isolated points of $X_h$ is countable, the scattered space $X_h$ is separable. 

It remains to prove that the space $X_h$ is Lindel\"of. Let $\U$ be any open cover of $X_h$. For every $x\in X$, find a set $U_x\in\U$ containing $x$ and an open neighborhood $V_x$ of $x$ in $K$ such that $$B_h(x;X\cap V_x)\defeq\{x\}\cup\{y\in X\cap V_x:h(y)<h(x)\}\subseteq U_x.$$ Observe that  $K\setminus \bigcup_{x\in Q}V_x$ is a compact subset of $\w^\w$ and hence it is contained in the lower set ${\downarrow}g$ of some function $g\in\w^\w$. Since the set $\{f_\alpha\}_{\alpha\in\kappa}$ is $\K$-Lusin, the set $\Omega=\{\alpha\in \kappa:f_\alpha\le^* g\}$ is countable. For every $\alpha\in\Omega$ the ordinal $\lambda_{\alpha+1}$ has countable cofinality and hence $\lambda_{\alpha+1}=\sup_n\lambda_{\alpha,n}$ for some strictly increasing sequence of ordinals $(\lambda_{\alpha,n})_{n\in\w}$ with $\lambda_{\alpha,0}=\lambda_\alpha$. Since the compact metrizable space $K$ is hereditarily Lindel\"of, for every $n\in\w$ there exists a countable set $Z_{\alpha,n}\subseteq D_{\lambda_{\alpha,n}}$ such that $D_{\lambda_{\alpha,n}}\subseteq\bigcup_{x\in Z_{\alpha,n}}V_x$.
Consider the countable subfamily
$$\U'\defeq\{U_x:x\in Q\cup \textstyle \bigcup_{\alpha\in\Omega} \bigcup_{n\in\w}Z_{\alpha,n}\}$$of the cover $\U$. We claim that the set $$\textstyle S\defeq X\setminus\bigcup\U'$$ is countable.

It follows from $Q=h^{-1}[{\uparrow}\lim\lambda]$ that  
$\bigcup_{x\in Q}V_x=\bigcup_{x\in Q}B_{h}(x;V_x)\subseteq\bigcup_{x\in Q}U_x$ and hence $S\subseteq K\setminus\bigcup_{x\in Q}V_x\subseteq{\downarrow}g$
and $$S\subseteq X\cap{\downarrow}g\subseteq\textstyle\bigcup_{\alpha\in\Omega}\big({\Uparrow} f_\alpha\cap{\downarrow}g\big).$$

Assuming that the set $S$ is uncountable, we can find an ordinal $\delta\in\Omega$ such that $S\cap{\Uparrow}f_\delta$ is uncountable. We can assume that $\delta$ is the smallest ordinal with this property. Then for every $\gamma\in\Omega\cap\delta$ the set $S\cap{\Uparrow}f_\gamma$ is countable and hence the set $$S\cap{\downarrow}g\cap{\Uparrow}f_\delta\setminus\textstyle\bigcup_{\gamma\in\Omega\cap\delta}{\Uparrow} f_\gamma$$ is uncountable. It follows that 
$$\textstyle G\defeq{\downarrow}g\cap{\Uparrow}f_\delta\setminus\bigcup_{\gamma\in\Omega\cap\delta}{\Uparrow}f_\gamma={\downarrow}g\cap{\Uparrow}f_\delta\setminus\bigcup_{\gamma\in\delta}{\Uparrow} f_\gamma$$ is a Borel and hence analytic subset of $\w^\w$ such that $G\cap S$ is uncountable. We claim that for every $n\in\w$ the closed subset 
$$\textstyle\Gamma_n\defeq G\setminus \bigcup_{x\in Z_{\delta,n}}V_x$$ 
of the analytic space $G$ is countable. In the opposite case $\Gamma_n$, contains an uncountable compact subset (see \cite[29.1]{Ke}) and hence $\Gamma_n$ has nonempty intersection with the Bernstein set $B_{\lambda_{\delta,n}}$, which is not possible as $B_{\lambda_{\delta,n}}\cap \Gamma_n\subseteq B_{\lambda_{\delta,n}}\cap G\subseteq D_{\lambda_{\delta,n}}\subseteq\bigcup_{x\in Z_{\delta,n}}V_x\subseteq K\setminus \Gamma_n$. 
This contradiction shows that for every $n\in\w$ the set $\Gamma_n$ is countable and so is the set $\Gamma\defeq\bigcup_{n\in\w}\Gamma_n$. Observe that 
$$\textstyle S\cap G\subseteq X\cap{\Uparrow} f_\delta\setminus \bigcup_{\gamma<\delta}{\Uparrow} f_\gamma=\bigcup_{\lambda_\delta\le\alpha<\lambda_{\delta+1}}D_\alpha$$
and $h[S\cap G]\subseteq [\lambda_\delta,\lambda_{\delta+1})$. Pick any element $y\in S\cap G\setminus \Gamma$ and find $n\in\w$ such that $h(y)<\lambda_{\delta,n}$. 
It follows from $y\notin \Gamma$ that $y\in G\setminus \Gamma_n\subseteq \bigcup_{x\in Z_{\delta,n}}V_x$ and hence $y\in V_x$ for some $x\in Z_{\delta,n}\subseteq D_{\lambda_{\alpha,n}}$. Taking into account that $h(y)<\lambda_{\delta,n}=h(x)$, we conclude that $y\in B_h(x;V_x)\subseteq U_x\in \U'$, which contradicts $y\in S$. This contradiction completes the proof of the Lindel\"of property of the space $X_h$.
\end{proof}

\begin{theorem}\label{t:construct4} If there exists a $\K$-Lusin set of uncountable cardinality $\kappa$, then for every ordinal $\lambda$ with $2\le \lambda<\kappa^+$ and $\cf(\lambda)\le\w$ there exists a separable special FLS space $X_h$ of cardinality $\kappa$ and scattered height $\lambda$. If $\kappa=|\lambda|$, then the space $X_h$ is hereditarily separable.
\end{theorem}

\begin{proof}  {Assume that there exists a $\K$-Lusin set of cardinality $\kappa>\w$. By (the proof of the equivalence $(1)\Leftrightarrow(2)$ in) Lemma~\ref{l:bQK}, there exists a $\IQ$-concentrated set $Z\subseteq\IR$} of cardinality $\kappa$. Since subsets of $\IQ$-concentrated sets are $\IQ$-concentrated, we can assume that $0\in Z$ and $x-y\notin\IQ+\IZ\sqrt{2}$ for every distinct elements $x,y\in Z$. Let $\lambda$ be any ordinal such that $2\le\lambda<\kappa^+$ and $\cf(\lambda)\le\w$. Separately we shall consider the cases of successor and limit ordinal $\lambda$. 
\smallskip

First, assume that the ordinal $\lambda$ is successor. Since $\lambda<\kappa^+$, we can find a cover $(Z_\alpha)_{\alpha\in\lambda}$ of $Z$ by pairwise disjoint nonempty subsets such that $Z_{\max\lambda}=\{0\}$. If $\kappa=|\lambda|$, then we can (and do) additionally assume that each set $Z_\alpha$ is countable. Let $X\defeq Z+\IQ$ and $h:X\to\lambda$ be the function defined by $h^{-1}\{\alpha\}=Z_\alpha+\IQ$ for $\alpha\in\lambda$. Since each set $h^{-1}\{\alpha\}$ is dense in $X$, the function $h$ is a height function. By Proposition~\ref{p:Xh1}, the scattered modification $X_h$ of $X$ is first-countable and has scattered height $\lambda$. It is clear that $|X|=\kappa$. For the largest element $\alpha$ of the successor ordinal $\lambda$, the subspace $h^{-1}[{\uparrow}\alpha]=h^{-1}\{\alpha\}=Z_{\alpha}+\IQ=\{0\}+\IQ$ of $X_h$ is countable and hence Lindel\"of. Since the set $Z$ is $\IQ$-concentrated, the space $X=Z+\IQ\subseteq\IR$ is concentrated at $h^{-1}[{\uparrow}\alpha]=\IQ$. By Lemma~\ref{l:XhL1}, the space $X_h$ is Lindel\"of.
\smallskip

Next, assume that the ordinal $\lambda$ is limit. Since $\cf(\lambda)\le\w$, there exists a strictly increasing sequence of ordinals $(\lambda_n)_{n\in\w}$ such that $\lambda_0=0$ and $\sup_n\lambda_n=\lambda$.  Let $\{Z_\alpha\}_{\alpha\in\lambda}$ be a cover of $Z$ by pairwise disjoint nonempty sets such that $Z_0=\{0\}$ and for every $n\in\w$ the set $Z_{{\lambda}_n}$ is a singleton.  If $\kappa=|\lambda|$, then we can (and will) additionally assume that each set $Z_\alpha$ is countable. 
Let 
$$\textstyle X\defeq (\IQ+\w\sqrt{2})\cup\bigcup_{n\in\w}\bigcup_{\lambda_n<\alpha< \lambda_{n+1}}(Z_\alpha+(n+1)\sqrt{2}+\IQ)$$
and $h:X\to\lambda$ be the function such that for every $n\in\w$ and ordinal $\alpha\in\lambda$ with $\lambda_n<\alpha<\lambda_{n+1}$ we have $h^{-1}\{\lambda_n\}=\IQ+n\sqrt{2}$ and $h^{-1}\{\alpha\}=Z_{\alpha}+(n+1)\sqrt{2}+\IQ$. Since for every $\alpha\in\lambda$ the fiber $h^{-1}\{\alpha\}$ is dense in $X$, the function $h$ is a height function. By Proposition~\ref{p:Xh1}, the scattered modification $X_h$ of $X$ is a first-countable scattered space of scattered height $\lambda$. To show that $X_h$ is Lindel\"of, we shall apply Lemma~\ref{l:XhL2}. First observe that the subspace $h^{-1}[{\uparrow}\lambda]$ of $X_h$ is empty and hence Lindel\"of. For every $n\in\w$ the subspace $h^{-1}\{\lambda_{n+1}\}=\IQ+(n+1)\sqrt{2}$ of $X$ is countable and hence Lindel\"of. Also the space $$h^{-1}[{\uparrow}\lambda_n,{\downarrow}\lambda_{n+1}]\subseteq (\IQ+n\sqrt{2})\cup (Z+\IQ+(n+1)\sqrt{2})$$ is concentrated at $\IQ+(n+1)\sqrt{2}$ (because $Z$ is concentrated at $\IQ$). By Lemma~\ref{l:XhL2}, the space $X_h$ is Lindel\"of. Therefore, $X_h$ is a special FLS space of cardinality $\kappa$ and scattered height $\lambda$.

If $\kappa=|\lambda|$, then the countability of the sets $Z_\alpha$ implies the countability of fibers of the height function $h:X\to \lambda$. Now Proposition~\ref{p:Xh2}(5) ensures that the FLS space $X_h$ is hereditarily separable.
\end{proof}





\section{Properties of FLS spaces}

In this section we explore properties of FLS spaces and show that the constructions of FLS spaces presented in the preceding section are close to being the best possible. First we derive simple upper bounds on the cardinality and scattered height of FLS spaces. 

\begin{proposition}\label{p:ImpossibleRanks} If $X$ is an FLS space, then  $|X|\le\mathfrak c$, $\hbar[X]<\mathfrak c^+$ and $\cf(\hbar[X])\le\w$.
\end{proposition}

\begin{proof}
By Arhangel'ski\u{\i}'s inequality \cite{Arhangelskii}, every first-countable Lindel\"of space has cardinality at most $\continuum$, which implies $|\hbar[X]|\le|X|\le\continuum$ and hence $\hbar[X]<\continuum^+$. The scattered height of a Lindel\"{o}f scattered space $X$ cannot have uncountable cofinality, because  $\set{X \setminus X^{(\a)}}{\a < \hbar[X]}$ is an open cover for $X$ with no subcovers of size $< \cf(\hbar[X])$.
\end{proof}

Proposition~\ref{p:ImpossibleRanks} shows that the FLS spaces constructed in Theorems~\ref{thm:Existence}, \ref{t:construct2}, \ref{t:construct3} have the maximum possible cardinality. 

\begin{question}\label{q:AllSizes?}
Is there an FLS space of cardinality $\k$ for every $\k \leq \continuum$?
\end{question}

Theorem~\ref{t:construct4} and Remark~\ref{rem:Qd} below imply that it is consistent with arbitrarily large values of $\continuum$ that the answer to Question~\ref{q:AllSizes?} is positive. Moreover, the results in this section show that if we restrict the question to FLS spaces of finite scattered height, then the answer is independent of \zfc. But we do not know whether a negative answer to Question~\ref{q:AllSizes?} is also consistent. 
The next theorem and its corollary provide a partial solution: Martin's Axiom implies that every submetrizable FLS space is either countable or has cardinality $\continuum$.

\begin{theorem}\label{t:countable} For an FLS space $X$ the following conditions are equivalent:
\begin{enumerate}
\item[\textup{(1)}] $X$ is countable;
\item[\textup{(2)}] $X$ is a $Q$-space;
\item[\textup{(3)}] $X$ is perfect; 
\item[\textup{(4)}] $X$ is hereditarily Lindel\"of;
\item[\textup{(5)}] $X$ is symmetrizable;
\item[\textup{(6)}] $X$ has countable separating weight and $|X|<\mathfrak q_1$.
\end{enumerate}
If $X$ is functionally Hausdorff, then the conditions \textup{(1)--(6)} are equivalent to:
\begin{enumerate}
\item[\textup{(7)}] $X$ is submetrizable and has cardinality $|X|<\mathfrak q_0$.
\end{enumerate}
\end{theorem}

\begin{proof} The implications $(1)\Ra(2)\Ra(3)$ are trivial.
\smallskip

$(3)\Ra(4)$ Assume that $X$ is perfect. To show that $X$ is hereditarily Lindel\"of, take any subspace $Y\subseteq X$ and any  cover $\U$ of $Y$ by open subsets of $X$. Since the space $X$ is perfect, the open set $\bigcup\U$ can be written as the union $\bigcup_{n\in\w}F_n$ of a sequence $(F_n)_{n\in\w}$ of closed subsets of $X$. For every $n\in\w$ the closed subset $F_n$ of the Lindel\"of space $X$ is Lindel\"of and hence $F_n\subseteq \bigcup\U_n$ for some countable subfamily $\U_n$ of $\U$. Then $\bigcup_{n\in\w}\U_n$ is a countable subcover of $Y$, witnessing that the space $Y$ is Lindel\"of and $X$ is hereditarily Lindel\"of.
\smallskip

$(4)\Ra(1)$ Assume that the space $X$ is hereditarily Lindel\"of. Then for every ordinal $\alpha\in\hbar[X]$ the discrete subspace $X^{(\alpha)}\setminus X^{(\alpha+1)}$ of $X$ is countable. By the Lindel\"of property of the subspace $X\setminus X^{(\w_1)}$, there exists a countable ordinal $\alpha$ such that $X\setminus X^{(\w_1)}=X\setminus X^{(\alpha)}$, which means that $X^{(\alpha)}=X^{(\w_1)}=\emptyset$ and hence the scattered height of $X$ is countable. Therefore, the space $X=\bigcup_{\alpha\in\hbar[X]}(X^{(\alpha)}\setminus X^{(\alpha+1)})$ is countable being a countable union of countable sets.
\smallskip

$(1)\Ra(5)$ Assume that $X$ is countable. Then $X$ admits a linear order $\le$ such that for every $x\in X$ the set $\{y\in X:y\le x\}$ is finite. For every $x\in X$ fix a neighborhood base $\{U_n(x)\}_{n\in\w}$ at $x$ such that $U_{n+1}(x)\subseteq U_n(x)\subseteq U_0(x)=X$ for all $n\in\w$. It is easy to see that the symmetric $d:X\times X\to[0,\infty)$ defined by $$d(x,y)\defeq\max\{2^{-n}:n\in\w\mbox{ and }\max\{x,y\}\in U_n(\min\{x,y\})\}$$generates the topology of $X$.
\smallskip

$(5)\Ra(3)$ Assume that $X$ is symmetrizable and let $d$ be a symmetric generating the topology of $X$. By the first-countability of $X$, for every $x\in X$ and $\e>0$ the ball $B_d(x;\e)=\{y\in X:d(x,y)<\e\}$ is a neighborhood of $x$, see \cite[9.6]{Grue}. Consequently, $x$ is contained in the interior $B_d(x;\e)^\circ$ of $B_d(x;\e)$ in $X$. This implies that every closed subset $F$ of $X$ is equal to the intersection $\bigcap_{n\in\w}\bigcup_{x\in F}B_d(x;2^{-n})^\circ$ and hence is of type $G_\delta$ in $X$, witnessing that $X$ is perfect.
\smallskip

$(1)\Ra(6)$ Assume that $X$ is countable and for every distinct points $x,y\in X$ choose an open set $U_{x,y}$ that contains $x$ but not $y$. Then $\{U_{x,y}:x,y\in X,\;x\ne y\}$ is a countable separating family witnessing that $X$ has countable separating weight. The strictly inequality $|X|<\mathfrak q_1$ is trivial.
\smallskip

$(6)\Ra(2)$ Assume that $X$ has countable separating weight and $|X|<\mathfrak q_1$. Let $\U$ be a countable family of open sets separating points of the space $X$ and let $\tau$ be a second-countable topology on $X$, generated by the subbase $\U$. The separating property of $\U$ implies that the second-countable topological space $X_\tau=(X,\tau)$ is $T_1$. By the definition of the cardinal $\mathfrak q_1$, the space $X_\tau$ is a $Q$-space and so is the space $X$.
\smallskip

Now assuming that $X$ is functionally Hausdorff, we shall prove that $(1)\Ra(7)\Ra(2)$.

$(1)\Ra(7)$:  For any pair $(a,b)$ in the set $P\defeq\{(x,y)\in X\times X:x\ne y\}$, choose a continuous function $f_{a,b}:X\to\IR$ such that $f_{a,b}(a)\ne f_{a,b}(b)$.
Observe that the diagonal map $f:X\to\IR^P$, $f:x\mapsto (f_{a,b}(x))_{(a,b)\in P}$ is injective. If $X$ is countable, then so is the set $P$ and the space $\IR^P$ is metrizable, witnessing that $X$ is submetrizable. The equality $|X|\le\w<\mathfrak q_0$ is trivial.
\smallskip

$(7)\Ra(2)$: Assume that $X$ is submetrizable and $|X|<\mathfrak q_0$. Find a  continuous bijective map $f:X\to Y$ onto a metrizable space $Y$. Since $X$ is Lindel\"of, the metrizable space $Y$ is Lindel\"of and separable. The metrizable separable space $f[X]$ has cardinality $<\mathfrak q_0$ and hence is a $Q$-space by the definition of the cardinal $\mathfrak q_0$. Then for every set $A\subseteq X$ it image $f[A]$ is of type $F_\sigma$ in $Y$. By the continuity of $f$, the preimage $f^{-1}[f[A]]=A$ of $f[A]$ is an $F_\sigma$-set in $X$, witnessing that $X$ is a $Q$-space.
\end{proof}

\begin{corollary}\label{p:Qsets} 
Every submetrizable FLS space of cardinality less than $\mathfrak{q}_0$ is countable. In particular, Martin's Axiom implies that every submetrizable FLS space is either countable or of cardinality $\continuum$.
\end{corollary}

Theorem~\ref{thm:Existence} provides many \zfc-examples of uncountable FLS spaces of infinite scattered height. The problem of the existence of an uncountable FLS of finite scattered height is more subtle. First we consider the question of what cardinalities such spaces can have. 

\begin{proposition}\label{p:AllSizes}
Let $\k$ be a cardinal, and suppose there is a cardinality-$\k$ FLS space with finite scattered height. Then there is a cardinality-$\lambda$ FLS space with finite scattered height for every $\lambda \leq \k$.
\end{proposition}

\begin{proof}
Suppose $X$ is an FLS space with $|X| = \k$ and $\hbar[X] < \w$. There is a largest integer $r < \hbar[X]$ such that $|X^{(r)}| \geq \lambda$. Let $Y$ be any subset of $X^{(r)}$ with $|Y| = \lambda$, and let $Z = \closure{Y}$. Then $Y \sub Z \sub Y \cup X^{(r+1)}$, and this implies $|Z| = \lambda$. And because $Z$ is a closed subspace of $X$, $Z$ is an FLS space with finite scattered height. 
\end{proof}

In light of {Proposition~\ref{p:AllSizes}}, the question of what are the possible cardinalities of finite-height FLS spaces really reduces to the question of what is the largest possible cardinality of a finite-height FLS space. We shall show that such largest possible cardinality is $\mathfrak d$.

We shall say that a subset $A$ of a topological space is {\em of type} $G_\kappa$ in $X$ where $\kappa$ is a cardinal if $A=\bigcap_{\alpha\in\kappa}U_\alpha$ for some family $\{U_\alpha\}_{\alpha\in\kappa}$ of open sets in $X$. 

\begin{lemma}\label{l:fsh} Let $X$ be a topological space of countable extent and countable pseudocharacter. Then every closed scattered subspace $Y\subseteq X$ of finite scattered height has cardinality $|Y|\le\mathfrak d$ and is a $G_{\mathfrak d}$-set in $X$. If $\mathfrak b>\w_1$, then every closed subspace $Y\subseteq X$ of finite scattered height is countable.
\end{lemma}

\begin{proof} The lemma will be proved by induction on the scattered height $\hbar[Y]$. If $\hbar[Y]=0$, then $Y=\emptyset$ and hence $|Y|=0\le\mathfrak d$ and $Y=\emptyset$ is an open (and hence $G_{\mathfrak d}$-set) in $X$.

Assume that for some positive integer number $k$ we have proved that every closed subspace of scattered height $<k$ in $X$ has cardinality $\le\mathfrak d$ and is a $G_{\mathfrak d}$-set in $X$. Let $Y\subseteq X$ be a closed subspace of scattered height $k$. It follows that $Y^{(k-1)}$ is a non-empty closed discrete subspace of $X$. Since $X$ has countable extent, the set $Y^{(k-1)}$ is countable and hence can be written as $Y^{(k-1)}=\{y_n:n\in\w\}$ for some sequence $(y_n)_{n\in\w}$. 

Since $X$ has countable pseudocharacter, for every $x\in X$ there exists a decreasing sequence of open sets $(U_n(x))_{n\in\w}$ such that $\bigcap_{n\in\w}U_n(x)=\{x\}$. For every $g\in\w^\w$, consider the open neighborhood $V_g\defeq\bigcup_{n\in\w}U_{g(n)}(y_n)$ of the set $Y^{(k-1)}$ in $X$.

 For every $x\in X\setminus Y^{(k-1)}$ consider the function $f_x\in\w^\w$ defined by $f_x(n)=\min\{m\in\w:x\notin U_m(y_n)\}$ and observe that $x\notin V_g$ for every function $g\in{\uparrow}f_x$.

By definition of the cardinal $\mathfrak d$, there exists a subset $\mathcal D\subseteq\w^\w$ of cardinality $|\mathcal D|=\mathfrak d$ such that $\w^\w=\bigcup_{g\in\mathcal D}{\downarrow}g$. For every function $g\in\mathcal D$, the subspace $$Y_g\defeq Y\setminus V_g$$ of $Y$ is closed and has scattered height $<k$. By the inductive hypothesis, $|Y_g|\le\mathfrak d$ and $Y_g$ is a $G_{\mathfrak d}$-set in $X$. Consequently, there exists a family $\U_g$ of open sets in $X$ such that $Y_g=\bigcap\U_g$ and $|\U_g|\le\mathfrak d$.  

Consider the family $$\V\defeq\big\{U\cup V_g:g\in\mathcal D,\;\;U\in\U_g\big\}$$of open neighborhoods of the set $Y$ in $X$ and observe that $|\V|\le\mathfrak d$. We claim that $Y=\bigcap\V$. Indeed, for any $x\in X\setminus Y$, we can find a function $g\in\mathcal D$ such that $f_x\le g$, which implies $x\notin V_g$. Since $x\notin Y_g=\bigcap\U_g$, there exists a set $U\in\U_f$ such that $x\notin U$. Then $x\notin U\cup V_g\in\V$ and hence $x\notin\bigcap\V$, witnessing that $Y$ is a $G_{\mathfrak d}$-set in $X$.

The dominating property of $\mathcal D$ ensures that $Y\setminus Y^{(k-1)}\subseteq\bigcup_{g\in\mathcal D}Y_g$ and hence
$$\textstyle |Y|=\big|Y^{(k-1)}\cup\bigcup_{g\in\mathcal D}Y_g\big|\le \mathfrak d.$$

Now assuming that $\mathfrak b>\w_1$, we shall prove that every closed subspace $Y$ of finite scattered height in $X$ is countable. Since $X$ has countable extent, the nonempty closed discrete subset $Y^{(\hbar[Y]-1)}$ of $Y$ and $X$ is countable. Let $l<\hbar[Y]$ be the smallest number such that the set $Y^{(l)}$ is countable. If $l=0$, then $Y=Y^{(0)}$ is countable and we are done. So, assume that $l>0$ and conclude that $Y^{(l-1)}$ is uncountable. Since $l<\hbar[Y]$, the countable set $Y^{(l)}$ is not empty and hence can be written as $Y^{(l)}=\{y_n:n\in\w\}$ for some sequence $(y_n)_{n\in\w}$. For every $y\in Y^{(l-1)}\setminus Y^{(l)}$ and every $n\in\w$, let $f_y(n)=\min\{k\in\w:y\notin U_k(y_n)\}$. Observe that for every $f\in\w^\w$ the set 
$$\textstyle \{y\in Y^{(l-1)}\setminus Y^{(l)}:f_y\le f\}=Y\setminus\bigcup_{n\in\w}U_{f(n)}(y_n)$$ 
is countable, being a closed discrete subspace of the space $X$. This implies that for every $f\in\w^\w$,  the set $\{y\in Y^{(l-1)}\setminus Y^{(l)}:f_y\le^* f\}$ is countable and hence any uncountable subset of $\{f_y:y\in Y^{(l-1)}\setminus Y^{(l)}\}$ is $\le^*$-unbounded in $\w^\w$ witnessing that $\mathfrak b=\w_1$. 
\end{proof}

\begin{corollary}\label{c:bd} 
Every FLS space of finite scattered height has cardinality $\leq\mathfrak d$. If $\mathfrak b>\w_1$, then every FLS space of finite scattered height is countable.
\end{corollary}

Next, we derive some restrictions on the cofinality of the limit part of the scattered height of an FLS space.

\begin{proposition}\label{p:cf-gL} If for an uncountable cardinal $\kappa$ there exists a  Lindel\"of scattered space $X$ of countable pseudocharacter with $\cf(\lim\hbar[X])=\kappa$ and $|X^{(\lim\hbar[X])}|\le\w$, then there exists a generalized $\K$-Lusin set of cardinality $\kappa$.
\end{proposition}

\begin{proof} Let $\lambda=\lim\hbar[X]$ and observe that $\lambda<\hbar[X]$ (as $\cf(\hbar[X])\le\w<\kappa=\cf(\lambda)$). Then the countable set $X^{(\lambda)}$ is not empty and hence $X^{(\lambda)}=\{y_n\}_{n\in\w}$ for some sequence $(y_n)_{n\in\w}$. Since the space $X$ has countable pseudocharacter, for every point $x\in X$ there exists a decreasing sequence $(U_n(x))_{n\in\w}$ of open sets in $X$ such that $\bigcap_{n\in\w}U_n(x)=\{x\}$. Let $(\lambda_\alpha)_{\alpha < \kappa}$ be an increasing sequence of ordinals with $\sup_{\alpha\in\kappa}\lambda_\alpha=\lambda$. For each $\alpha < \kappa$, choose a point $x_\alpha$ in the set $X^{(\lambda_\alpha+1)}\setminus X^{(\lambda_\alpha)}$ and define a function $f_\alpha \in \w^\w$ by setting $f_\alpha(n) = \min\{m\in\w: x_\alpha \notin U_m(y_n)\}$ for each $n\in\w$.

We claim that the set $L\defeq\{f_\alpha:\alpha\in\kappa\}\subseteq\w^\w$ has cardinality $\kappa$ and is generalized $\K$-Lusin. Given any $g\in\w^\w$, consider the open neighborhood $V_g\defeq\bigcup_{n\in\w}U_{g(n)}(y_n)$ of the set $X^{(\lambda)}$ in $X$. Since the space $X$ is Lindel\"of, the open cover $V_g\cup\{X\setminus X^{(\alpha)}:\alpha\in\lambda\}$ has a countable subcover. Since $\cf(\lambda)=\kappa>\w$, there exists  an ordinal $\beta\in\lambda$ such that $X=V_g\cup (X\setminus X^{(\beta)})$ and thus $X^{(\beta)}\subseteq V_g$. Then for every ordinal $\alpha\in\kappa$ with $\hbar(x_\alpha)\ge \beta$ we have $x_\alpha\in X^{(\beta)}\subseteq V_g=\bigcup_{n\in\w}U_{g(n)}(y_n)$ and hence $x_\alpha\in U_{g(n)}(y_n)$ for some $n\in\w$. Now the definition of the function $f_\alpha$ ensures that $g(n)<f_\alpha(n)$ and hence $f_\alpha\not\le g$. The choice of the increasing sequence $(\lambda_\alpha)_{\alpha\in\kappa}$ with $\sup_{\alpha\in\kappa}\lambda_\alpha=\lambda$ ensures that the set $$\{\alpha\in\kappa:f_\alpha\le g\}\subseteq\{\alpha\in\kappa:\hbar(x_\alpha)<\beta\}=\{\alpha\in\kappa:\lambda_\alpha<\beta\}$$ has cardinality $<\kappa$. This implies that the function $f_*:\kappa\to \w^\w$, $f_*:\alpha\mapsto f_\alpha$, has fibers of cardinality $<\kappa$ and by the regularity of the cardinal $\kappa=\cf(\lambda)$, the set $L=f_*[\kappa]$ has cardinality $\kappa$. To see that $L$ is generalized $\K$-Lusin, take any compact subset $K\subseteq\w^\w$ and find a function $g\in\w^\w$ such that $K\subseteq {\downarrow}g$. As we already know,  the set $\{\alpha\in\kappa:f_\alpha\le g\}$ has cardinality $<\kappa$, which implies that the set $L\cap K\subseteq L\cap{\downarrow}g=\{f_\alpha:\alpha\in\kappa,\;f_\alpha\le g\}$ also has cardinality $<\kappa=|L|$ and this means that $L$ is generalized $\K$-Lusin.
\end{proof}

\begin{corollary} If $\mathfrak b>\w_1$ and $X$ is an FLS space with $\kappa=\cf(\lim\hbar[X])>\w$, then there exists a generalized $\K$-Lusin set of cardinality $\kappa$.
\end{corollary}

\begin{proof} Let $\lambda=\lim\hbar[X]$. Since $X^{(\lambda)}$ is an FLS space of finite scattered height, $X^{(\lambda)}$ is countable by Corollary~\ref{c:bd}. Now we can apply Proposition~\ref{p:cf-gL} and conclude that there exists a generalized $\K$-Lusin set of cardinality $\kappa$.
\end{proof}

\begin{question} Can the equality $|X^{(\lim\hbar[X])}|\le\w$ be removed from the requirements of Proposition~\ref{p:cf-gL}? {\rm Observe that for a Lindel\"of scattered space $X$, the set $X^{(\lim\hbar[X])}$ is countable if $\hbar[X]=\lim\hbar[X]+1$ or $\mathfrak b>\w_1$.}
\end{question}

\begin{proposition}\label{p:bd} If $X$ is a scattered Lindel\"of space $X$ of  countable pseudocharacter, then $$\cf(\lim\hbar[X])\in\{0,\w\}\cup [\mathfrak b,\mathfrak d].$$
\end{proposition}

\begin{proof} Let $\lambda=\lim\hbar[X]$. Assuming that $\cf(\lambda)>\w$, we should prove that $\mathfrak b\le\cf(\lambda)\le\mathfrak d$. By Proposition~\ref{p:ImpossibleRanks}, $\cf(\hbar[X])\le\w$ and hence $\lambda<\hbar[X]$. Then the subspace $X^{(\lambda)}$ of $X$ is nonempty, closed and has finite scattered height. By Lemma~\ref{l:fsh}, $X^{(\lambda)}$ is a $G_{\mathfrak d}$-set in the Lindel\"of space $X$. This implies that every open cover of $X\setminus X^{(\lambda)}$ has a subcover of cardinality $\le\mathfrak d$. Taking into account that the open cover $\{X\setminus X^{(\alpha)}:\alpha\in\lambda\}$ of $X\setminus X^{(\lambda)}$ does not have a subcover of cardinality less than $\cf(\lambda)$, we conclude that $\cf(\lambda)\le\mathfrak d$.

It remains to show that $\cf(\lambda)\ge\mathfrak b$. To derive a contradiction, assume that $\cf(\lambda)<\mathfrak b$. Then $\mathfrak b>\w_1$ and the set $X^{(\lambda)}$ is countable according to Lemma~\ref{l:fsh}. By Proposition~\ref{p:cf-gL}, there exists a generalized $\K$-Lusin set $L$ of uncountable regular cardinality $\cf(\lambda)<\mathfrak b$, which contradicts Lemma~\ref{l:Kbd}(1).
\end{proof} 

Lemma~\ref{l:fsh} implies that uncountable FLS spaces of finite scattered height can exist only under the set-theoretic assumption $\mathfrak b=\w_1$. Now we show that the equality $\mathfrak b=\w_1$ is also sufficient for the existence of uncountable FLS spaces with finite scattered height. 

\begin{theorem}\label{t:main} For an uncountable cardinal $\kappa$ the following conditions are equivalent:
\begin{enumerate}
\item[\textup{(1)}] for every ordinal $\lambda$ with $2\le\lambda<\kappa^+$ and 
$\cf(\lambda)\le\w$ there exists a hereditarily separable special FLS space of cardinality $\kappa$ and scattered height $\lambda$;
\item[\textup{(2)}] there exists an FLS space of cardinality $\kappa$ and scattered height $2$;
\item[\textup{(3)}] there exists a $\K$-Lusin subset of cardinality $\kappa$.
\end{enumerate}
\end{theorem}

\begin{proof} The implication $(3)\Ra(1)$ was proved in Theorem~\ref{t:construct4} and $(1)\Ra(2)$ is trivial.
\smallskip

$(2)\Ra(3)$: Assume that $X$ is an FLS space of cardinality $\kappa$ and scattered height 2. Since $\hbar[X]=2$, the set $X^{(1)}$ of non-isolated points of $X$ is a non-empty closed discrete subspace in $X$. Since $X$ is Lindel\"of, its closed discrete subspace $X^{(1)}$ is countable and  hence $X^{(1)}=\{x_n\}_{n\in\w}$ for some sequence $(x_n)_{n\in\w}$. For every $x\in X$ fix a decreasing sequence $(U_n(x))_{n\in\w}$ of open sets such that $\{x\}=\bigcap_{n\in\w}U_n(x)$. For every $x\in X\setminus X^{(1)}$ consider the function $f_x\in\w^\w$ defined by $f_x(n)=\min\{m\in\w:x\notin U_m(x_n)\}$ for every $n\in\w$. We claim that $L=\{f_x:x\in X\setminus X^{(1)}\}$ is a $\K$-Lusin set of cardinality $\kappa$. Observe that for every $g\in\w^\w$ the set 
$$\textstyle \{x\in X\setminus X^{(1)}:f_x\le g\}=X\setminus \bigcup_{n\in\w}U_{g(n)}(x_n)$$ 
is countable, being a closed  discrete subspace of the Lindel\"of space $X$. This implies that the function $f_*:X\setminus X^{(1)}\to \w^\w$, $f_*:x\mapsto f_x$, has countable fibers and hence $|L|=|f_*[X\setminus X^{(1)}]|=|X\setminus X^{(1)}|=|X|=\kappa$. To see that the set $L$ is $\K$-Lusin, take any compact set $K\subseteq\w^\w$ and find a function $g\in\w^\w$ such that $K\subseteq{\downarrow}g$. Since the set $X\setminus\bigcup_{n\in\w}U_{g(n)}(x_n)$ is countable, so is the set $\textstyle L\cap K\subseteq\{f_x:x\in X\setminus X^{(1)},\;f_x\le g\}=\{f_x:x\in X\setminus\bigcup_{n\in\w}U_{g(n)}(x_n)\}.$
\end{proof}

\begin{corollary}\label{c:b} The following conditions are equivalent:
\begin{enumerate}
\item[\textup{(1)}] $\mathfrak b=\w_1$;
\item[\textup{(2)}] there exists an uncountable $\IQ$-concentrated set in the real line;
\item[\textup{(3)}] there exists an uncountable $\K$-Lusin set;
\item[\textup{(4)}] for every ordinal $\lambda$ with $2\le\lambda<\w_2$, 
$\cf(\lambda)\le\w$ and $\cf(\lim\lambda)\le\w_1$ there exists a uncountable FLS-space of scattered height $\lambda$;
\item[\textup{(5)}] there exists an uncountable FLS-space of finite scattered height;
\item[\textup{(6)}] there exists an FLS-space $X$ with $\cf(\lim \hbar[X])=\w_1$;
\item[\textup{(7)}] there exists a semiregular uncountable FLS space of scattered height $2$.
\end{enumerate}
\end{corollary}

\begin{proof} The equivalence of the first three conditions was proved in Lemma~\ref{l:bQK}. The implication $(3)\Ra(4)$ follows from Theorem~\ref{t:main} and the implications $(4)\Ra(5,6)$ are trivial. The implication $(5)\Ra(1)$ follows from Corollary~\ref{c:bd}, and $(6)\Ra(1)$ from Proposition~\ref{p:bd}.
The implication $(7)\Ra(5)$ is trivial and $(4)\Ra(7)$ is established in the following proposition.
\end{proof}

\begin{proposition} Every FLS space $X$ of scattered height $2$ has a point-countable base and is a closed subspace of a semiregular FLS space $Y$ with $|Y|=|X|$ and with scattered height $2$.
\end{proposition}

\begin{proof} Since $X$ is Lindel\"of and has scattered height 2, the closed subset $X^{(1)}$ of $X$ is discrete and hence countable. Since $X$ is first-countable, every point $x\in X$ has a countable neighborhood base $\{U_n(x)\}_{n\in\w}$. Then $$\{\{x\}:x\in X\setminus X^{(1)}\}\cup\{U_n(x):n\in\w,\;x\in X^{(1)}\}$$is a point-countable base for $X$.

Next, we construct an embedding of $X$ into a semiregular FLS space. Let $Y$ be any set containing $X$ so that $Y\setminus X$ is countable. Write $Y\setminus X$ as the union $\bigcup_{n\in\w}F_n$ of an increasing sequence $(F_n)_{n\in\w}$ of finite sets. Let $(S_x)_{x\in X^{(1)}}$ be a disjoint family of infinite subsets of $Y\setminus X$. On the set $Y$, consider the topology generated by the point-countable base
$$\mathcal B=\{\{y\}:y\in Y\setminus X^{(1)}\}\cup\{U_n(x)\cup S_x\setminus F_n:x\in X^{(1)},\;n\in\w\}.$$
It is clear that $Y$ is an FLS space of scattered height 2 and $X$ is a closed subspace in $Y$.

It remains to show that the space $Y$ is semiregular. For this it suffices to check that every basic set $U\in\mathcal B$ is regular open. This is clear if $U=\{x\}$ for some $x\in Y\setminus X^{(1)}$. So, assume that $U=U_n(x)\cup S_x\setminus F_n$  for some $x\in X^{(1)}$ and $n\in\w$. Assuming that $U$ is not regular open, find an interior point $y\in\overline{U}\setminus U$ of the closure $\overline{U}$ of $U$ in $Y$. Since each point of the set $Y\setminus X^{(1)}$ is isolated, the point $y$ belongs to the set $X^{(1)}=Y^{(1)}$.  Since $y$ is an interior point of $\overline{U}$, there exists $m\in\w$ such that $U_m(y)\cup S_y\setminus F_m\subseteq \overline{U}$. On the other hand, the definition of the topology of the space $Y$ ensures that  $\overline{U}\subseteq X\cup S_x$ and thus $\overline{U}$ is disjoint with the infinite set $S_y\setminus F_m\subseteq\overline{U}$. This contradiction shows that each set $U\in\mathcal B$ is regular open and the space $Y$ is semiregular.
\end{proof}

Now let us ask some questions motivated by Theorem~\ref{t:main} and its Corollary~\ref{c:b}. By Corollary~\ref{c:b} the existence of an uncountable FLS space of finite scattered space is equivalent to the existence of an uncountable $\K$-Lusin set. On the other hand, by Theorem~\ref{t:main} the existence of a $\K$-Lusin set of cardinality $\kappa$ is equivalent to the existence of a size-$\kappa$ FLS of scattered height 2.
We do not know whether the ``scattered height $2$'' in Theorem~\ref{t:main}(3) can be weakened to ``finite scattered height''.

\begin{question}\label{q:Rank2}
Suppose there is a size-$\k$ FLS space with finite scattered height. Does this imply the equivalent conditions in Theorem~\ref{t:main}?
\end{question}

\begin{question}
In particular, is it consistent that there is an FLS space with height $3$ and cardinality $\geq\! \w_2$, but that every FLS space with height $2$ has cardinality $\le\!\w_1$?
\end{question}

As a partial answer to Question~\ref{q:Rank2}, we end this section by proving that the existence of a size-$\kappa$ FLS space of finite scattered height implies the existence of a generalized $\K_\sigma$-Lusin set of cardinality $\kappa$.

\begin{proposition}\label{p:GeneralizedLusin} If for an uncountable cardinal $\kappa$ there exists an FLS space of cardinality $\kappa$ and finite scattered height, then there exists a generalized $\mathcal K_\sigma$-Lusin set of cardinality $\k$.
\end{proposition}

\begin{proof} 
Suppose there is a cardinality-$\k$ FLS space of finite scattered height, and let $r$ be the minimal scattered height of such a space. 
Let $X$ be an FLS space with $\hbar[X] = r$ and $|X| = \k$. 
Because $X^{(r-1)}$ is a closed discrete subset of the Lindel\"of space $X$, $X^{(r-1)}$ is countable and hence can be written as $X^{(r-1)} = \set{x_n}{n \in \w}$ for some sequence $(x_n)_{n\in\w}$.

Since $X$ has countable pseudocharacter, for every $x\in X$ there exists a decreasing sequence $(U_n(x))_{n\in\w}$ of open sets such that $\{x\}=\bigcap_{n\in\w}U_n(x)$. For every function $g\in\w^\w$ consider the open neighborhood $V_g\defeq\bigcup_{n\in\w}U_{g(n)}(x_n)$ of $X^{(r-1)}$. Observe that for every countable set $G\subseteq\w^\w$, the $F_\sigma$-subspace $X\setminus \bigcap_{g\in G}V_g$ of $X$ is an FLS space of scattered height $<r$. The minimality of $r$ ensures that $|X\setminus \bigcap_{g\in G}V_g|<\kappa$.

For each $x\in X\setminus X^{(r-1)}$ consider the function $f_x\in\w^\w$ defined by $f_x(n)=\min\{m\in\w:x\not\in U_m(x_n)\}$. Choose a family $(\hat f_x)_{x\in X\setminus X^{(r-1)}}$ of pairwise distinct functions such that $f_x\le \hat f_x$ for every $x\in X\setminus X^{(r-1)}$. Then the set $L\defeq\{\hat f_x:x\in X\setminus X^{(r-1)}\}$ has cardinality $\kappa$. We claim that $L$ is generalized $\K_\sigma$-Lusin. Given any $\sigma$-compact set $K\subseteq\w^\w$, find a countable set $G\subseteq\w^\w$ such that $K\subseteq\bigcup_{g\in G}{\downarrow}g$ and observe that 
$$\textstyle L\cap K\subseteq \bigcup_{g\in G}\{\hat f_x:x\in X\setminus X^{(r-1)},\;\hat f_x\le g\}\subseteq \bigcup_{g\in G}\{\hat f_x:x\in X\setminus X^{(r-1)},\;f_x\le g\}$$
 and hence $|L\cap K|\le \big|\bigcup_{g\in G}\{x\in X\setminus X^{(r-1)}:f_x\le g\}\big|=\big|X\setminus \bigcap_{g\in G}V_g\big|<\kappa$.
\end{proof}

\section{Characterizing possible cardinalities and heights of FLS spaces}

In this section we present several results describing the possible cardinalities and scattered heights of FLS spaces under various set-theoretic assumptions. (Some of this is simply summarizing what we have proved in the preceding sections.) 
We also include some consistency results along these lines.

\begin{corollary}\label{t:size-c} Suppose there exist $\K$-Lusin sets of all cardinalities $<\mathfrak c$. Then for any cardinal $\k$, there is an FLS space of cardinality $\k$ if and only if $\kappa\le\mathfrak c$.
\end{corollary}

\begin{proof} The forward implication follows from Proposition~\ref{p:ImpossibleRanks}, and the reverse implication follows from Theorems~\ref{thm:Existence} and \ref{t:construct4}.
\end{proof}

The hypothesis of this corollary is consistent by Remark~\ref{rem:Qd} below. Recall from Corollary~\ref{p:Qsets} that the assumption $\mathfrak q_0=\mathfrak c$ (a consequence of $\ma$) has a totally different effect:

\begin{corollary} If $\mathfrak q_0 = \mathfrak c$, then there is a submetrizable FLS space of {infinite} cardinality $\k$
if and only if $\kappa\in\{\w,\mathfrak c\}$.
\end{corollary}

Next, we characterize possible cardinalities of FLS spaces with finite scattered height.

\begin{corollary}\label{c:size-k} Assume that for some uncountable cardinal $\lambda$ there exists a $\K$-Lusin set of cardinality $\lambda$ and there are no generalized $\K$-Lusin sets of cardinality $\lambda^+$. Then for any cardinal $\kappa$, there is a finite-height FLS space with cardinality $\k$ if and only if $\k\leq\lambda$.
\end{corollary}

\begin{proof} 
The reverse implication follows from Theorem~\ref{t:construct4}.

For the forward implication, assume, aiming for a contradiction, that $\kappa$ is equal to the cardinality of some FLS space of finite scattered height but $\kappa>\lambda$. Then $\lambda^+\le\kappa$ and by Proposition~\ref{p:AllSizes}, there exists an FLS space of cardinality $\lambda^+$ and finite scattered height. By Proposition~\ref{p:GeneralizedLusin}, there exists a generalized $\K$-Lusin set of cardinality $\lambda^+$, contradicting our assumptions concerning $\lambda$.
\end{proof}

\begin{corollary}\label{c:size-b} Suppose $\mathfrak b = {\w}_1$ and there are no generalized $\K$-Lusin sets of cardinality ${\w}_2$. Then for any cardinal $\kappa$, there is a finite-height FLS space with cardinality $\k$ if and only if $\kappa\le {\w}_1$.
\end{corollary}

\begin{proof} This follows from Corollary~\ref{c:size-k} and Lemma~\ref{l:bQK}.
\end{proof}

We note that the hypotheses of Corollary~\ref{c:size-b} are consistent by Theorem~\ref{thm:DualCohen} below.

\begin{corollary}\label{c:size-d1} Suppose there is a $\K$-Lusin set of cardinality $\mathfrak d$. Then for any cardinal $\kappa$, there is a finite-height FLS space with cardinality $\k$ if and only if $\kappa\le\mathfrak d$.
\end{corollary}

\begin{proof} This follows from Corollary~\ref{c:bd}, Corollary~\ref{c:size-k}, and Lemma~\ref{l:Kbd}.
\end{proof}

The hypothesis of Corollary~\ref{c:size-d1} is consistent by Remark~\ref{rem:Qd} below.

\begin{corollary}\label{c:size-d} If $\mathfrak b > {\w}_1$, then every finite-height FLS space is countable.
\end{corollary}

\begin{proof} This follows from Corollary~\ref{c:bd}.
\end{proof}

Finally, we characterize possible scattered heights of FLS spaces.

\begin{corollary}\label{c:height} If there exists a $\K$-Lusin set of cardinality $\mathfrak d$, then for any ordinal $\lambda$ the following conditions are equivalent:
\begin{enumerate}
\item[\textup{(1)}] $\lambda$ is equal to the scattered height of some uncountable FLS space;
\item[\textup{(2)}] $2\le\lambda<\mathfrak c^+$, $\cf(\lambda)\le\w$ and $\cf(\lim \lambda)\le \mathfrak d$.
\end{enumerate}
\end{corollary}

\begin{proof} The implication $(1)\Ra(2)$ follows from Proposition~\ref{p:ImpossibleRanks} and \ref{p:bd}, and $(2)\Ra(1)$ follows from Theorems~\ref{t:construct3} (for infinite $\lambda$) and \ref{t:construct4} (for finite $\lambda$).
\end{proof} 

\begin{corollary}\label{c:height1} If $\mathfrak d = {\w}_1$, then for any ordinal $\lambda$ the following conditions are equivalent:
\begin{enumerate}
\item[\textup{(1)}] $\lambda$ is equal to the scattered height of some uncountable FLS space;
\item[\textup{(2)}] $2\le\lambda<\mathfrak c^+$, $\cf(\lambda)\le\w$ and $\cf(\lim \lambda)\le \w_1$.
\end{enumerate}
\end{corollary}

\begin{proof} This corollary follows from Corollary~\ref{c:height} and Lemma~\ref{l:genL}.
\end{proof}

\begin{corollary}\label{c:height2} If $\w_1<\mathfrak b=\mathfrak d=\cov{M}$, then for any  ordinal $\lambda$ the following conditions are equivalent:
\begin{enumerate}
\item[\textup{(1)}] $\lambda$ is equal to the scattered height of some uncountable FLS space;
\item[\textup{(2)}] $\w\le\lambda<\mathfrak c^+$, $\cf(\lambda)\le\w$ and $\cf(\lim \lambda)\in\{\w,{\cf(\mathfrak d)}\}$.
\end{enumerate}
\end{corollary}

\begin{proof} $(1)\Ra(2)$ Assume that $\lambda=\hbar[X]$ for some uncountable FLS space $X$. Since $\w_1<\cov{M}=\mathfrak b=\mathfrak d$, Corollary~\ref{c:bd} implies that $\w\le \lambda$. By Propositions~\ref{p:ImpossibleRanks} and \ref{p:bd}, $\lambda<\mathfrak c^+$, $\cf(\lambda)\le\w$ and $\cf(\lim\lambda)\in\{\w\}\cup[\mathfrak b,\mathfrak d]$. {By \cite[2.4]{Blass}, $\mathfrak b\le\cf(\mathfrak d)\le\mathfrak d$ and hence $\mathfrak b=\mathfrak d$ implies $\cf(\lim\lambda)\in\{\w\}\cup[\mathfrak b,\mathfrak d]=\{\w,{\cf(\mathfrak d)}\}$.}
\smallskip

The implication $(2)\Ra(1)$ follows from Theorems~\ref{thm:Existence} and \ref{t:construct2}.
\end{proof}

\begin{remark}\label{rem:Qd} By \cite[8.2.6]{BJ}, for any cardinal $\kappa$ of uncountable cofinality, adding $\kappa$ many Cohen reals to a ground model of \ch produces a model of \zfc in which $\w^\w$ contains a Lusin (and hence $\K$-Lusin) set of cardinality $\kappa=\mathfrak c$. In this model, $\mathfrak b=\w_1$ and $\mathfrak d=\mathfrak c$. In Theorem~\ref{thm:CohenModel} we shall construct a model of \zfc in which $\w_1=\mathfrak b<\mathfrak d<\mathfrak c$ and there exists a $\K$-Lusin set of cardinality $\mathfrak d$. In Theorem~\ref{thm:DualCohen} we shall construct a model of \zfc in which $\w_1=\mathfrak b<\w_2<\mathfrak d=\mathfrak c$ and there are no generalized $\K$-Lusin sets of cardinality $\w_2$, which means that the cardinal $\kappa$ in Corollary~\ref{c:size-k} can be equal to $\w_1$, which is strictly smaller than $\mathfrak d$.
\end{remark}
 
 Theorem~\ref{thm:Existence}  and Corollaries~\ref{c:height}--\ref{c:height2} motivate the following intriguing question.

\begin{question}\label{q:possibleheights1}
Is it consistent that $\cf(\lim\hbar[X])\le\w$ for every FLS space $X$? 
\end{question}

\begin{question}\label{q:possibleheights}
If $\mathfrak b = \mathfrak d = \k$, is there an FLS space with scattered height $\k+1$?
\end{question}

In combination with these results, the following two theorems show:
\begin{itemize}
\item[$\bullet$] It is consistent with arbitrary values of $\mathfrak d$ and $\mathfrak c$ that there is a finite-height FLS space of cardinality $\k$ for every $\k \leq \mathfrak{d}$.
\item[$\bullet$] It is consistent with arbitrary values of $\mathfrak{d}$ that there are finite-height FLS spaces of cardinality $\aleph_1$, but no larger.
\end{itemize}

\begin{theorem}\label{thm:CohenModel}
Let $\lambda$ and $\mu$ be cardinals with uncountable cofinality such that $\continuum \leq \lambda \leq \mu$. There is a ccc forcing extension in which $\dom = \lambda$, $\continuum = \mu$, and there exists a $\K$-Lusin set of cardinality $\mathfrak d$.
\end{theorem}

\begin{proof}
Let $\CC$ denote the standard poset for adding $\lambda$ Cohen reals. In $V^{\CC}$, let $\BB$ be the standard poset for adding $\mu$ random reals. We claim that $V^{\CC * \dot \BB}$ satisfies the conclusions of the theorem.

Standard arguments (which can be found, for example, in \cite{BJ}) show that $\dom = \lambda$ and $\continuum = \mu$ in $V^{\CC*\dot \BB}$. Thus,  we will be done if we can show that 
there exists in $V^{\CC*\dot\BB}$ a size-$\lambda$ set $B \sub \w^\w$ such that every uncountable subset of $B$ is $\leq^*$-unbounded.

In $V^{\CC}$ there is a $\K$-Lusin set $B\subseteq\w^\w$ of size $\lambda$ in $\w^\w$. (Any uncountable subset of the mutually generic Cohen reals added by $\CC$ is Lusin; see \cite[Section 8.2]{BJ} for details.) 
We claim that the set $B$ remains $\K$-Lusin in $V^{\CC*\dot\BB}$.
This follows from the fact that the poset $\BB$ has the $\w^\w$-bounding property: a subset of $\w^\w$ in $V^{\PP}$ is $\leq^*$-bounded in $V^{\PP}$ if and only if it is $\leq^*$-bounded in $(V^{\PP})^\BB = V^{\PP*\dot\BB}$.
\end{proof}

\begin{theorem}\label{thm:DualCohen}
Let $\lambda \geq \max\{\continuum,\aleph_3\}$ be a regular uncountable cardinal. There is a ccc forcing extension in which $\aleph_1 = \bdd < \dom = \continuum =  \lambda$, there is no generalized $\K$-Lusin sets of cardinality $\aleph_2$ and every uncountable FLS space of finite scattered height has cardinality $\aleph_1$. 
\end{theorem}

\begin{proof}
Let $\PP$ be a ccc poset that forces Martin's Axiom \ma plus $\continuum = \lambda$, and let $\CC = \mathrm{Fn}(\w_1,2)$ be the standard poset for adding $\aleph_1$ Cohen reals. We claim that $V^{\PP * \CC}$ satisfies the conclusions of the theorem.

Standard arguments (which can be found, for example, in \cite[Chapter 3]{BJ}) show that $\bdd = \aleph_1$ and $\continuum = \dom = \lambda$ in $V^{\PP*\CC}$. 
Thus, applying Corollaries~\ref{c:b} and \ref{c:size-k}, we will be done if we can show that in $V^{\PP*\CC}$ there are no generalized $\mathcal K$-Lusin sets of size $\aleph_2$.

Working in $V^{\PP*\CC}$, suppose $X \sub \w^\w$ with $|X| = \aleph_2$. We aim to show $X$ is not a generalized $\mathcal K$-Lusin set.

For each $\a < \w_1$, let $\CC_\a = \mathrm{Fn}(\a,2)$ denote the restriction of $\CC$ to $\a$. 
For every $x \in X$, because $\CC$ has the ccc there is some $\a < \w_1$ such that $x \in V^{\PP*\CC_\a}$. 
By the Pigeonhole Principle, there is some particular $\a < \w_1$ such that $\big|X \cap V^{\PP*\CC_\a}\big| = \aleph_2$. 
Let $Y = X \cap V^{\PP*\CC_\a}$.

Because $\PP*\CC_\a$ has the ccc, there is in $V^{\PP*\CC_\a}$ some $Z \sub \w^\w$ such that $Z \supseteq Y$ and $|Z| = |Y| = \aleph_2$.
By a theorem of Roitman \cite{Roitman}, \masc holds in $V^{\PP*\CC_\a}$. (Specifically, Roitman proved that adding a single Cohen real to a model of \ma results in a model of \masc. This applies to $V^{\PP*\CC_\a}$ because $\CC_\a$ is equivalent to the forcing to add a single Cohen real.) But $\continuum = \lambda$ already in $V^{\PP}$, which implies $\continuum = \lambda$ in $V^{\PP*\CC_\a}$, and \masc implies $\bdd = \continuum$. Thus $\bdd > \aleph_2$ in $V^{\PP*\CC_\a}$. In particular, there is in $V^{\PP*\CC_\a}$ some $g \in \w^\w$ such that $f \leq^* g$ for every $f \in Z$.

In $V^{\PP*\CC}$, it is still true that $f \leq^* g$ for every $f \in Z$, and consequently $f \leq^* g$ for every $f \in Y$. By Lemma~\ref{lem:RemoveTheStar}, there is some $h \in \w^\w$ such that $\card{Y \cap \set{f \in \w^\w}{f \leq h}} = \aleph_2$. But $\set{f \in \w^\w}{f \leq h}$ is a compact subset of $\w^\w$, so this shows $Y$ is not a generalized $\mathcal K$-Lusin set.
\end{proof}




\begin{thebibliography}{99}

\bibitem{Arhangelskii} A. V. Arhangel'ski\u{\i}, ``On the cardinality of bicompacta satisfying the first axiom of countability,'' \emph{Soviet Mathematics Doklady} \textbf{10} (1969), 951--955.

\bibitem{BB} T.~Banakh, L.~Bazylevych, `{`$Q$-spaces, perfect spaces and related cardinal characteristics of the continuum'', \emph{preprint} ({\tt arxiv.org/abs/2206.01667})}.

\bibitem{BMZ} T. Banakh, M. Machura, and L. Zdomskyy, ``On critical cardinalities related to $Q$-sets,'' \emph{Mathematical Bulletin of Taras Shevchenko Scientific Society} \textbf{11} (2014), 21--32.

\bibitem{BJ} T. Bartoszynski and H. Judah, \emph{Set Theory: On the Structure of the Real Line}, A K Peters (1995).


\bibitem{Blass} A.~Blass, {``Combinatorial cardinal characteristics of the continuum'', \emph{Handbook of set theory}.} Vols. 1, 2, 3, 395--489, Springer, Dordrecht, 2010.


\bibitem{Engelking} R. Engelking, \emph{General Topology}. Sigma Series in Pure Mathematics, 6, Heldermann, Berlin (revised edition), 1989.

\bibitem{FM} W. G. Fleissner and A. W. Miller, ``On $Q$-sets,'' \emph{Proceedings of the American Mathematical Society} \textbf{78}:2 (1980), 280--284. 

\bibitem{Gewand} M. E. Gewand, ``The Lindel\"{o}f degree of scattered spaces and their products,'' \emph{Journal of the Australian Mathematical Society (Series A)} \textbf{37} (1984), 98--105.

\bibitem{Grue} G.~Gruenhage, {``Generalized metric spaces'', \emph{Handbook of set-theoretic topology},} 423--501, North-Holland, Amsterdam, 1984.

\bibitem{Hodel} R.~Hodel,{``Cardinal functions. I'', \emph{Handbook of set-theoretic topology},} 1--61, North-Holland, Amsterdam, 1984.

\bibitem{Ke} A.~Kechris, {\em Classical Descriptive Set Theory}, Springer, 1995.

\bibitem{Miller} A. W. Miller, ``Special subsets of the real line,'' in \emph{Handbook of {set-theoretic t}opology}, K. Kunen and J. E. Vaughan (eds.), North-Holland (1984), 201--234.

\bibitem{MM} M. N.~Mukherjee, D.~Mandal, {``Concerning nearly metrizable spaces'', \emph{Applied General Topology}} {\bf 14}:2 (2013), 135--145.

\bibitem{Oxtoby} J. Oxtoby, \emph{Measure and Category. A survey of the analogies between topological and measure spaces.} Second edition. Graduate Texts in Mathematics, 2. Springer-Verlag, New York-Berlin, 1980.

\bibitem{Roitman} J. Roitman, ``Adding a random or a Cohen real: topological consequences and the effect on Martin's axiom,'' \emph{Fundamenta Mathematicae} \textbf{103} (1979), 47--60.


\bibitem{SS} L.~Steen, J.~Seebach, Jr. {\em Counterexamples in topology}, Holt, Rinehart and Winston, Inc., New York-Montreal, Que.-London 1970.


\end{thebibliography}
\end{document}